\documentclass{article}
\usepackage{amsmath}
\usepackage{amssymb}
\usepackage{amsthm}
\usepackage{amscd}
\begin{document}
\theoremstyle{plain}
\newtheorem{thm}{Theorem}[section]
\newtheorem{lem}[thm]{Lemma}
\newtheorem{prop}[thm]{Proposition}
\newtheorem{cor}[thm]{Corollary}
\theoremstyle{definition}
\newtheorem{dfn}[thm]{Definition}
\newtheorem{ex}[thm]{Example}
\def\e#1\e{\begin{equation}#1\end{equation}}
\def\ea#1\ea{\begin{align}#1\end{align}}
\def\eq#1{{\rm(\ref{#1})}}
\def\dim{\mathop{\rm dim}}
\def\Re{\mathop{\rm Re}}
\def\Im{\mathop{\rm Im}}
\def\id{\mathop{\rm id}}
\def\ind{\mathop{\rm ind}}
\def\Hol{{\textstyle\mathop{\rm Hol}}}
\def\Ker{\mathop{\rm Ker}}
\def\Coker{\mathop{\rm Coker}}
\def\Vol{\mathop{\rm Vol}}
\def\rank{\mathop{\rm rank}}
\def\diam{\mathop{\rm diam}}
\def\sign{\mathop{\rm sign}}
\def\vol{\mathop{\rm vol}}
\def\O{\mathbin{\rm O}}
\def\SO{\mathbin{\rm SO}}
\def\GL{\mathbin{\rm GL}}
\def\U{\mathbin{\rm U}}
\def\SL{\mathop{\rm SL}}
\def\SU{\mathop{\rm SU}}
\def\Sp{\mathop{\rm Sp}}
\def\Spin{\mathop{\rm Spin}}
\def\sech{{\textstyle\mathop{\rm sech}}}
\def\ge{\geqslant} 
\def\le{\leqslant} 
\def\R{\mathbin{\mathbb R}}
\def\H{\mathbin{\mathbb H}}
\def\Z{\mathbin{\mathbb Z}}
\def\C{\mathbin{\mathbb C}}
\def\CP{\mathbb{CP}}
\def\al{\alpha}
\def\be{\beta}
\def\la{\lambda}
\def\ga{\gamma}
\def\de{\delta}
\def\io{\iota}
\def\ep{\epsilon}
\def\ka{\kappa}
\def\th{\theta}
\def\vp{\varphi}
\def\si{\sigma}
\def\ze{\zeta}
\def\De{\Delta}
\def\La{\Lambda}
\def\Om{\Omega}
\def\Ga{\Gamma}
\def\Si{\Sigma}
\def\om{\omega}
\def\d{{\rm d}}
\def\pd{\partial}
\def\db{{\bar\partial}}
\def\ts{\textstyle}
\def\sst{\scriptscriptstyle}
\def\w{\wedge}
\def\lt{\ltimes}
\def\sm{\setminus}
\def\op{\oplus}
\def\ot{\otimes}
\def\bigot{\bigotimes}
\def\iy{\infty}
\def\ra{\rightarrow}
\def\longra{\longrightarrow}
\def\hookra{\hookrightarrow}
\def\t{\times}
\def\ha{{\textstyle\frac{1}{2}}}
\def\ti{\tilde}
\def\ovB{\,\overline{\!B}}
\def\ms#1{\vert#1\vert^2}
\def\bms#1{\bigl\vert#1\bigr\vert^2}
\def\md#1{\vert #1 \vert}
\def\bmd#1{\big\vert #1 \big\vert}
\def\nm#1{\Vert #1 \Vert}
\def\cnm#1#2{\Vert #1 \Vert_{C^{#2}}} 
\def\lnm#1#2{\Vert #1 \Vert_{L^{#2}}} 
\def\bnm#1{\bigl\Vert #1 \bigr\Vert}
\def\bcnm#1#2{\bigl\Vert #1 \bigr\Vert_{C^{#2}}} 
\def\blnm#1#2{\bigl\Vert #1 \bigr\Vert_{L^{#2}}} 
\def\an#1{\langle#1\rangle}
\def\ban#1{\bigl\langle#1\bigr\rangle}
\title{The exceptional holonomy groups \\ and calibrated geometry}
\author{Dominic Joyce \\ Lincoln College, Oxford, OX1 3DR}
\date{}
\maketitle

\section{Introduction}
\label{ec1}

In the theory of Riemannian holonomy groups, perhaps the most
mysterious are the two exceptional cases, the holonomy group
$G_2$ in 7 dimensions and the holonomy group $\Spin(7)$ in 8 
dimensions. This is a survey paper on the exceptional holonomy
groups, in two parts. Part I collects together useful facts
about $G_2$ and $\Spin(7)$ in \S\ref{ec2}, and explains
constructions of compact 7-manifolds with holonomy $G_2$ in
\S\ref{ec3}, and of compact 8-manifolds with holonomy $\Spin(7)$
in~\S\ref{ec4}.

Part II discusses the {\it calibrated submanifolds} of
manifolds of exceptional holonomy, namely {\it associative
$3$-folds} and {\it coassociative $4$-folds} in $G_2$-manifolds,
and {\it Cayley $4$-folds} in $\Spin(7)$-manifolds. We introduce
calibrations in \S\ref{ec5}, defining the three geometries and
giving examples. Finally, \S\ref{ec6} explains their {\it
deformation theory}.

Sections \ref{ec3} and \ref{ec4} describe my own work, for
which the main reference is my book \cite{Joyc5}. Part II
describes work by other people, principally the very important
papers by Harvey and Lawson \cite{HaLa} and McLean \cite{McLe},
but also more recent developments.

This paper was written to accompany lectures at the $11^{\rm th}$
G\"okova Geometry and Topology Conference in May 2004, sponsored
by TUBITAK. In keeping with the theme of the conference, I have
focussed mostly on $G_2$, at the expense of $\Spin(7)$. The paper
is based in part on the books \cite{Joyc5} and \cite[Part I]{GHJ},
and the survey paper~\cite{Joyc8}.
\medskip

\noindent{\it Acknowledgements.} I would like to thank the
conference organizers Turgut Onder and Selman Akbulut for
their hospitality. Many people have helped me develop my
ideas on exceptional holonomy and calibrated geometry;
I would particularly like to thank Simon Salamon and
Robert Bryant.
\vfil\eject

\centerline{\LARGE\scshape Part I. Exceptional Holonomy}
\medskip

\section{Introduction to $G_2$ and $\Spin(7)$}
\label{ec2}

We introduce the notion of {\it Riemannian holonomy groups},
and their classification by Berger. Then we give short
descriptions of the holonomy groups $G_2$, $\Spin(7)$ and
$\SU(m)$, and the relations between them. All the results
below can be found in my book~\cite{Joyc5}.

\subsection{Riemannian holonomy groups}
\label{ec21}

Let $M$ be a connected $n$-dimensional manifold, $g$ a Riemannian
metric on $M$, and $\nabla$ the Levi-Civita connection of $g$. Let
$x,y$ be points in $M$ joined by a smooth path $\ga$. Then {\it
parallel transport\/} along $\ga$ using $\nabla$ defines an isometry
between the tangent spaces $T_xM$, $T_yM$ at $x$ and~$y$.

\begin{dfn} The {\it holonomy group} $\Hol(g)$ of $g$ is the group
of isometries of $T_xM$ generated by parallel transport around
piecewise-smooth closed loops based at $x$ in $M$. We consider
$\Hol(g)$ to be a subgroup of $\O(n)$, defined up to conjugation
by elements of $\O(n)$. Then $\Hol(g)$ is independent of the base
point $x$ in~$M$.
\label{ec2def1}
\end{dfn}

Let $\nabla$ be the {\it Levi-Civita connection} of $g$. A tensor
$S$ on $M$ is {\it constant\/} if $\nabla S=0$. An important
property of $\Hol(g)$ is that it {\it determines the constant
tensors on}~$M$.

\begin{thm} Let\/ $(M,g)$ be a Riemannian manifold, and\/ $\nabla$
the Levi-Civita connection of\/ $g$. Fix a base point\/ $x\in M$,
so that\/ $\Hol(g)$ acts on $T_xM$, and so on the tensor powers
$\bigot^kT_xM\ot\bigot^lT_x^*M$. Suppose $S\in C^\iy\bigl(
\bigot^kTM\ot\bigot^lT^*M\bigr)$ is a constant tensor. Then
$S\vert_x$ is fixed by the action of\/ $\Hol(g)$. Conversely,
if\/ $S\vert_x\in\bigot^kT_xM\ot\bigot^lT_x^*M$ is fixed by
$\Hol(g)$, it extends to a unique constant tensor~$S\in
C^\iy\bigl(\bigot^kTM\ot\bigot^lT^*M\bigr)$.
\label{ec2thm1}
\end{thm}

The main idea in the proof is that if $S$ is a constant tensor and
$\ga:[0,1]\ra M$ is a path from $x$ to $y$, then $P_\ga(S\vert_x)=
S\vert_y$, where $P_\ga$ is the {\it parallel transport map} along
$\ga$. Thus, constant tensors are invariant under parallel transport.
In particular, they are invariant under parallel transport around
closed loops based at $x$, that is, under elements of~$\Hol(g)$.

The classification of holonomy groups was achieved by Berger 
\cite{Berg} in~1955.

\begin{thm} Let\/ $M$ be a simply-connected, $n$-dimensional 
manifold, and\/ $g$ an irreducible, nonsymmetric Riemannian 
metric on $M$. Then either
\begin{itemize}
\setlength{\parsep}{0pt}
\setlength{\itemsep}{0pt}
\item[{\rm(i)}] $\Hol(g)=\SO(n)$,
\item[{\rm(ii)}] $n=2m$ and\/ $\Hol(g)=\SU(m)$ or\/ $\U(m)$, 
\item[{\rm(iii)}] $n=4m$ and\/ $\Hol(g)=\Sp(m)$ or\/ $\Sp(m)\Sp(1)$, 
\item[{\rm(iv)}] $n=7$ and\/ $\Hol(g)=G_2$, or
\item[{\rm(v)}] $n=8$ and\/ $\Hol(g)=\Spin(7)$.
\end{itemize}
\label{ec2thm2}
\end{thm}

Here are some brief remarks about each group on Berger's list.

\begin{itemize}
\setlength{\labelwidth}{10pt}
\setlength{\parsep}{0pt}
\setlength{\itemsep}{0pt}
\item[(i)] $\SO(n)$ is the holonomy group of generic Riemannian metrics. 
\item[(ii)] Riemannian metrics $g$ with $\Hol(g)\subseteq\U(m)$ are 
called {\it K\"ahler metrics}. K\"ahler metrics are a natural class of 
metrics on complex manifolds, and generic K\"ahler metrics on a 
given complex manifold have holonomy $\U(m)$. 

Metrics $g$ with $\Hol(g)=\SU(m)$ are called
{\it Calabi--Yau metrics}. Since $\SU(m)$ is a subgroup of $\U(m)$, 
all Calabi--Yau metrics are K\"ahler. If $g$ is K\"ahler and $M$
is simply-connected, then $\Hol(g)\subseteq\SU(m)$ if and only 
if $g$ is Ricci-flat. Thus Calabi--Yau metrics are locally more
or less the same as Ricci-flat K\"ahler metrics.
\item[(iii)] metrics $g$ with $\Hol(g)=\Sp(m)$ are called
{\it hyperk\"ahler}. As $\Sp(m)\subseteq\SU(2m)\subset\U(2m)$,
hyperk\"ahler metrics are Ricci-flat and K\"ahler.

Metrics $g$ with holonomy group $\Sp(m)\Sp(1)$ for $m\ge 2$
are called {\it quaternionic K\"ahler}. (Note that quaternionic K\"ahler 
metrics are not in fact K\"ahler.) They are Einstein, but not Ricci-flat. 
\item[(iv),(v)] $G_2$ and $\Spin(7)$ are the exceptional cases,
so they are called the {\it exceptional holonomy groups}. Metrics
with these holonomy groups are Ricci-flat.
\end{itemize}

The groups can be understood in terms of the four {\it division 
algebras}: the {\it real numbers} $\R$, the {\it complex numbers} 
$\C$, the {\it quaternions} $\H$, and the {\it octonions} or 
{\it Cayley numbers} $\mathbb O$. 
\begin{itemize}
\item $\SO(n)$ is a group of automorphisms of $\R^n$. 
\item $\U(m)$ and $\SU(m)$ are groups of automorphisms of $\C^m$ 
\item $\Sp(m)$ and $\Sp(m)\Sp(1)$ are automorphism groups of $\H^m$. 
\item $G_2$ is the automorphism group of $\Im\mathbb{O}\cong\R^7$. 
$\Spin(7)$ is a group of automorphisms of ${\mathbb O}\cong\R^8$,
preserving part of the structure on~$\mathbb O$.
\end{itemize}

For some time after Berger's classification, the exceptional
holonomy groups remained a mystery. In 1987, Bryant \cite{Brya2}
used the theory of exterior differential systems to show that
locally there exist many metrics with these holonomy groups,
and gave some explicit, incomplete examples. Then in 1989, Bryant 
and Salamon \cite{BrSa} found explicit, {\it complete} metrics 
with holonomy $G_2$ and $\Spin(7)$ on noncompact manifolds.

In 1994-5 the author constructed the first examples of metrics
with holonomy $G_2$ and $\Spin(7)$ on {\it compact} manifolds
\cite{Joyc1,Joyc2,Joyc3}. These, and the more complicated
constructions developed later by the author \cite{Joyc4,Joyc5}
and by Kovalev \cite{Kova}, are the subject of Part~I.

\subsection{The holonomy group $G_2$}
\label{ec22}

Let $(x_1,\dots,x_7)$ be coordinates on $\R^7$. Write
$\d{\bf x}_{ij\ldots l}$ for the exterior form $\d x_i\w\d x_j
\w\cdots\w\d x_l$ on $\R^7$. Define a metric $g_0$, a 3-form $\vp_0$
and a 4-form $*\vp_0$ on $\R^7$ by~$g_0=\d x_1^2+\cdots+\d x_7^2$,
\e
\begin{split}
\vp_0&=\d{\bf x}_{123}+\d{\bf x}_{145}+\d{\bf x}_{167}
+\d{\bf x}_{246}-\d{\bf x}_{257}-\d{\bf x}_{347}-\d{\bf x}_{356}
\;\>\text{and}\\
*\vp_0&=\d{\bf x}_{4567}+\d{\bf x}_{2367}+\d{\bf x}_{2345}
+\d{\bf x}_{1357}-\d{\bf x}_{1346}-\d{\bf x}_{1256}-\d{\bf x}_{1247}.
\end{split}
\label{ec2eq1}
\e
The subgroup of $GL(7,\R)$ preserving $\vp_0$ is the {\it exceptional
Lie group} $G_2$. It also preserves $g_0,*\vp_0$ and the orientation
on $\R^7$. It is a compact, semisimple, 14-dimensional
Lie group, a subgroup of~$\SO(7)$.

A {\it $G_2$-structure} on a 7-manifold $M$ is a principal 
subbundle of the frame bundle of $M$, with structure group $G_2$. 
Each $G_2$-structure gives rise to a 3-form $\vp$ and a metric $g$ 
on $M$, such that every tangent space of $M$ admits an isomorphism 
with $\R^7$ identifying $\vp$ and $g$ with $\vp_0$ and $g_0$ 
respectively. By an abuse of notation, we will refer to $(\vp,g)$
as a $G_2$-structure.

\begin{prop} Let\/ $M$ be a $7$-manifold and\/ $(\vp,g)$ a
$G_2$-structure on $M$. Then the following are equivalent:
\begin{itemize}
\setlength{\parsep}{0pt}
\setlength{\itemsep}{0pt}
\item[{\rm(i)}] $\Hol(g)\subseteq G_2$, and\/ $\vp$ is the
induced\/ $3$-form,
\item[{\rm(ii)}] $\nabla\vp=0$ on $M$, where $\nabla$ is the
Levi-Civita connection of $g$, and
\item[{\rm(iii)}] $\d\vp=\d^*\vp=0$ on $M$.
\end{itemize}
\label{ec2prop1}
\end{prop}

Note that $\Hol(g)\subseteq G_2$ if and only if $\nabla\vp=0$
follows from Theorem \ref{ec2thm1}. We call $\nabla\vp$ the
{\it torsion} of the $G_2$-structure $(\vp,g)$, and when
$\nabla\vp=0$ the $G_2$-structure is {\it torsion-free}. A
triple $(M,\vp,g)$ is called a $G_2$-{\it manifold} if $M$ is a
7-manifold and $(\vp,g)$ a torsion-free $G_2$-structure on $M$.
If $g$ has holonomy $\Hol(g)\subseteq G_2$, then $g$ is Ricci-flat.

\begin{thm} Let\/ $M$ be a compact\/ $7$-manifold, and suppose 
that $(\vp,g)$ is a torsion-free $G_2$-structure on $M$. 
Then $\Hol(g)=G_2$ if and only if\/ $\pi_1(M)$ is finite. In 
this case the moduli space of metrics with holonomy $G_2$ on $M$, 
up to diffeomorphisms isotopic to the identity, is a smooth 
manifold of dimension~$b^3(M)$.
\label{ec2thm3}
\end{thm}

\subsection{The holonomy group $\Spin(7)$}
\label{ec23}

Let $\R^8$ have coordinates $(x_1,\dots,x_8)$. Define a 4-form
$\Om_0$ on $\R^8$ by
\e
\begin{split}
\Om_0=\,&\d{\bf x}_{1234}+\d{\bf x}_{1256}+\d{\bf x}_{1278}
+\d{\bf x}_{1357}-\d{\bf x}_{1368}-\d{\bf x}_{1458}-\d{\bf x}_{1467}\\
-&\d{\bf x}_{2358}-\d{\bf x}_{2367}-\d{\bf x}_{2457}
+\d{\bf x}_{2468}+\d{\bf x}_{3456}+\d{\bf x}_{3478}+\d{\bf x}_{5678}.
\end{split}
\label{ec2eq2}
\e
The subgroup of $\GL(8,\R)$ preserving $\Om_0$ is the holonomy group
$\Spin(7)$. It also preserves the orientation on $\R^8$ and the
Euclidean metric $g_0=\d x_1^2+\cdots+\d x_8^2$. It is a compact,
semisimple, 21-dimensional Lie group, a subgroup of~$\SO(8)$.

A $\Spin(7)$-structure on an 8-manifold $M$ gives rise to a 4-form 
$\Om$ and a metric $g$ on $M$, such that each tangent space of 
$M$ admits an isomorphism with $\R^8$ identifying $\Om$ and $g$ 
with $\Om_0$ and $g_0$ respectively. By an abuse of notation we 
will refer to the pair $(\Om,g)$ as a $\Spin(7)$-structure.

\begin{prop} Let\/ $M$ be an $8$-manifold and\/ $(\Om,g)$ a
$\Spin(7)$-structure on $M$. Then the following are equivalent:
\begin{itemize}
\setlength{\parsep}{0pt}
\setlength{\itemsep}{0pt}
\item[{\rm(i)}] $\Hol(g)\subseteq\Spin(7)$, and\/ $\Om$ is the 
induced\/ $4$-form,
\item[{\rm(ii)}] $\nabla\Om=0$ on $M$, where $\nabla$ is the
Levi-Civita connection of\/ $g$, and
\item[{\rm(iii)}] $\d\Om=0$ on $M$.
\end{itemize}
\label{ec2prop2}
\end{prop}

We call $\nabla\Om$ the {\it torsion} of the $\Spin(7)$-structure 
$(\Om,g)$, and $(\Om,g)$ {\it torsion-free} if $\nabla\Om=0$. 
A triple $(M,\Om,g)$ is called a $\Spin(7)$-{\it manifold} if $M$
is an 8-manifold and $(\Om,g)$ a torsion-free $\Spin(7)$-structure
on $M$. If $g$ has holonomy $\Hol(g)\subseteq\Spin(7)$, then $g$
is Ricci-flat.

Here is a result on {\it compact}\/ 8-manifolds with 
holonomy~$\Spin(7)$.

\begin{thm} Let\/ $(M,\Om,g)$ be a compact\/ $\Spin(7)$-manifold.
Then $\Hol(g)=\Spin(7)$ if and only if\/ $M$ is simply-connected,
and\/ $b^3(M)+b^4_+(M)=b^2(M)+2b^4_-(M)+25$. In this case the 
moduli space of metrics with holonomy $\Spin(7)$ on $M$, up to 
diffeomorphisms isotopic to the identity, is a smooth manifold of 
dimension\/~$1+b^4_-(M)$.
\label{ec2thm4}
\end{thm}

\subsection{The holonomy groups $\SU(m)$}
\label{ec24}

Let $\C^m\cong\R^{2m}$ have complex coordinates $(z_1,\ldots,z_m)$,
and define the metric $g_0$, K\"ahler form $\om_0$ and complex volume
form $\th_0$ on $\C^m$ by
\e
\begin{split}
g_0=\ms{\d z_1}+\cdots+\ms{\d z_m},\quad
\om_0&=\frac{i}{2}(\d z_1\w\d\bar z_1+\cdots+\d z_m\w\d\bar z_m),\\
\text{and}\quad\th_0&=\d z_1\w\cdots\w\d z_m.
\end{split}
\label{ec2eq3}
\e
The subgroup of $\GL(2m,\R)$ preserving $g_0,\om_0$ and $\th_0$
is the special unitary group $\SU(m)$. Manifolds with holonomy
$\SU(m)$ are called {\it Calabi--Yau manifolds}.

Calabi--Yau manifolds are automatically Ricci-flat and K\"ahler,
with trivial canonical bundle. Conversely, any Ricci-flat K\"ahler
manifold $(M,J,g)$ with trivial canonical bundle has $\Hol(g)
\subseteq\SU(m)$. By Yau's proof of the Calabi Conjecture
\cite{Yau}, we have:

\begin{thm} Let\/ $(M,J)$ be a compact complex $m$-manifold
admitting K\"ahler metrics, with trivial canonical bundle. Then
there is a unique Ricci-flat K\"ahler metric $g$ in each K\"ahler
class on $M$, and\/~$\Hol(g)\subseteq\SU(m)$.
\label{ec2thm5}
\end{thm}

Using this and complex algebraic geometry one can construct many
examples of compact Calabi--Yau manifolds. The theorem also applies in
the orbifold category, yielding examples of {\it Calabi--Yau orbifolds}.

\subsection{Relations between $G_2$, $\Spin(7)$ and $\SU(m)$}
\label{ec25}

Here are the inclusions between the holonomy groups $\SU(m),G_2$
and $\Spin(7)$:
\begin{equation*}
\begin{CD}
\SU(2) @>>> \SU(3) @>>> G_2 \\ @VVV @VVV @VVV \\
\SU(2)\t\SU(2) @>>> \SU(4) @>>> \Spin(7).
\end{CD}
\end{equation*}
We shall illustrate what we mean by this using the inclusion
$\SU(3)\hookra G_2$. As $\SU(3)$ acts on $\C^3$, it also acts
on $\R\op\C^3\cong\R^7$, taking the $\SU(3)$-action on $\R$ to
be trivial. Thus we embed $\SU(3)$ as a subgroup of $\GL(7,\R)$.
It turns out that $\SU(3)$ is a subgroup of the subgroup $G_2$
of $\GL(7,\R)$ defined in~\S\ref{ec22}.

Here is a way to see this in terms of differential forms. Identify
$\R\op\C^3$ with $\R^7$ in the obvious way in coordinates, so that
$\bigl(x_1,(x_2+ix_3,x_4+ix_5,x_6+ix_7)\bigr)$ in $\R\op\C^3$ is
identified with $(x_1,\ldots,x_7)$ in $\R^7$. Then $\vp_0=\d x_1
\w\om_0+\Re\th_0$, where $\vp_0$ is defined in \eq{ec2eq1} and
$\om_0,\th_0$ in \eq{ec2eq3}. Since $\SU(3)$ preserves $\om_0$
and $\th_0$, the action of $\SU(3)$ on $\R^7$ preserves $\vp_0$,
and so~$\SU(3)\subset G_2$.

It follows that if $(M,J,h)$ is Calabi--Yau 3-fold, then $\R\t M$ and
${\mathcal S}^1\t M$ have torsion-free $G_2$-structures, that is, are
$G_2$-manifolds.

\begin{prop} Let\/ $(M,J,h)$ be a Calabi--Yau $3$-fold, with K\"ahler
form $\om$ and complex volume form $\th$. Let\/ $x$ be a coordinate
on $\R$ or ${\mathcal S}^1$. Define a metric $g=\d x^2+h$ and a $3$-form
$\vp=\d x\w\om+\Re\th$ on $\R\t M$ or ${\mathcal S}^1\t M$. Then
$(\vp,g)$ is a torsion-free $G_2$-structure on $\R\t M$ or
${\mathcal S}^1\t M$, and\/~$*\vp=\ha\om\w\om-\d x\w\Im\th$.
\label{ec2prop3}
\end{prop}

Similarly, the inclusions $\SU(2)\hookra G_2$ and
$\SU(4)\hookra\Spin(7)$ give:

\begin{prop} Let\/ $(M,J,h)$ be a Calabi--Yau $2$-fold, with K\"ahler
form $\om$ and complex volume form $\th$. Let\/ $(x_1,x_2,x_3)$ be
coordinates on $\R^3$ or $T^3$. Define a metric $g=\d x_1^2+\d x_2^2+
\d x_3^2+h$ and a $3$-form $\vp$ on $\R^3\t M$ or $T^3\t M$ by
\e
\vp=\d x_1\w\d x_2\w\d x_3+\d x_1\w\om+\d x_2\w\Re\th-\d x_3\w\Im\th.
\label{ec2eq4}
\e
Then $(\vp,g)$ is a torsion-free $G_2$-structure on $\R^3\t M$ or
$T^3\t M$, and
\e
*\vp=\ha\om\w\om+\d x_2\w\d x_3\w\om-\d x_1\w\d x_3\w\Re\th
-\d x_1\w\d x_2\w\Im\th.
\label{ec2eq5}
\e
\label{ec2prop4}
\end{prop}

\begin{prop} Let\/ $(M,J,g)$ be a Calabi--Yau $4$-fold, with K\"ahler
form $\om$ and complex volume form $\th$. Define a $4$-form $\Om$ on
$M$ by $\Om=\ha\om\w\om+\Re\th$. Then $(\Om,g)$ is a
torsion-free $\Spin(7)$-structure on~$M$.
\label{ec2prop5}
\end{prop}

\section{Constructing $G_2$-manifolds from orbifolds $T^7/\Ga$\!\!}
\label{ec3}

We now explain the method used in \cite{Joyc1,Joyc2} and
\cite[\S 11--\S 12]{Joyc5} to construct examples of compact
7-manifolds with holonomy $G_2$. It is based on the {\it Kummer
construction} for Calabi--Yau metrics on the $K3$ surface, and
may be divided into four steps.

\begin{list}{}{\setlength{\leftmargin}{35pt}
\setlength{\labelwidth}{35pt}}
\item[Step 1.] Let $T^7$ be the 7-torus and $(\vp_0,g_0)$ a flat
$G_2$-structure on $T^7$. Choose a finite group $\Ga$ of isometries
of $T^7$ preserving $(\vp_0,g_0)$. Then the quotient $T^7/\Ga$ is
a singular, compact 7-manifold, an {\it orbifold}.

\item[Step 2.] For certain special groups $\Ga$ there is a 
method to resolve the singularities of $T^7/\Ga$ in a natural 
way, using complex geometry. We get a nonsingular, compact 
7-manifold $M$, together with a map $\pi:M\ra T^7/\Ga$, 
the resolving map.

\item[Step 3.] On $M$, we explicitly write down a 1-parameter 
family of $G_2$-structures $(\vp_t,g_t)$ depending on $t\in(0,\ep)$.
They are not torsion-free, but have small torsion when $t$ is small.
As $t\ra 0$, the $G_2$-structure $(\vp_t,g_t)$ converges to the
singular $G_2$-structure~$\pi^*(\vp_0,g_0)$.

\item[Step 4.] We prove using analysis that for sufficiently
small $t$, the $G_2$-structure $(\vp_t,g_t)$ on $M$, with
small torsion, can be deformed to a $G_2$-structure
$(\ti\vp_t,\ti g_t)$, with zero torsion. Finally, we show
that $\ti g_t$ is a metric with holonomy $G_2$ on the compact
7-manifold~$M$.
\end{list}

We will now explain each step in greater detail.

\subsection{Step 1: Choosing an orbifold}
\label{ec31}

Let $(\vp_0,g_0)$ be the Euclidean $G_2$-structure on $\R^7$
defined in \S\ref{ec22}. Suppose $\La$ is a {\it lattice} in
$\R^7$, that is, a discrete additive subgroup isomorphic to
$\Z^7$. Then $\R^7/\La$ is the torus $T^7$, and $(\vp_0,g_0)$
pushes down to a torsion-free $G_2$-structure on $T^7$. We
must choose a finite group $\Ga$ acting on $T^7$ preserving
$(\vp_0,g_0)$. That is, the elements of $\Ga$ are the push-forwards
to $T^7/\La$ of affine transformations of $\R^7$ which fix
$(\vp_0,g_0)$, and take $\La$ to itself under conjugation.

Here is an example of a suitable group $\Ga$, taken from
\cite[\S 12.2]{Joyc5}.

\begin{ex} Let $(x_1,\dots,x_7)$ be coordinates
on $T^7=\R^7/\Z^7$, where $x_i\in\R/\Z$. Let $(\vp_0,g_0)$ be
the flat $G_2$-structure on $T^7$ defined by \eq{ec2eq1}. Let
$\al,\be$ and $\ga$ be the involutions of $T^7$ defined by
\ea
\al:(x_1,\dots,x_7)&\mapsto(x_1,x_2,x_3,-x_4,-x_5,-x_6,-x_7),
\label{ec3eq1}\\
\be:(x_1,\dots,x_7)&\mapsto(x_1,-x_2,-x_3,x_4,x_5,\ha-x_6,-x_7),
\label{ec3eq2}\\
\ga:(x_1,\dots,x_7)&\mapsto\bigl(-x_1,x_2,-x_3,x_4,\ha-x_5,x_6,\ha-x_7).
\label{ec3eq3}
\ea
By inspection, $\al,\be$ and $\ga$ preserve $(\vp_0,g_0)$,
because of the careful choice of exactly which signs to change. 
Also, $\al^2=\be^2=\ga^2=1$, and $\al,\be$ and  $\ga$ commute.
Thus they generate a group $\Ga=\an{\al,\be,\ga}\cong\Z_2^3$ of
isometries of $T^7$ preserving the flat $G_2$-structure~$(\vp_0,g_0)$.
\label{ec3ex}
\end{ex}

Having chosen a lattice $\La$ and finite group $\Ga$, the quotient
$T^7/\Ga$ is an {\it orbifold}, a singular manifold with only quotient
singularities. The singularities of $T^7/\Ga$ come from the fixed points
of non-identity elements of $\Ga$. We now describe the singularities
in our example.

\begin{lem} In Example $\ref{ec3ex}$, $\be\ga,\ga\al,\al\be$ and\/
$\al\be\ga$ have no fixed points on\/ $T^7$. The fixed points of\/
$\al,\be,\ga$ are each\/ $16$ copies of\/ $T^3$. The singular set\/
$S$ of\/ $T^7/\Ga$ is a disjoint union of\/ $12$ copies of\/ $T^3$, 
$4$ copies from each of\/ $\al,\be,\ga$. Each component of\/ $S$
is a singularity modelled on that of\/~$T^3\t\C^2/\{\pm1\}$.
\end{lem}

The most important consideration in choosing $\Ga$ is that
we should be able to resolve the singularities of $T^7/\Ga$
within holonomy $G_2$. We will explain how to do this next.

\subsection{Step 2: Resolving the singularities}
\label{ec32}

Our goal is to resolve the singular set $S$ of $T^7/\Ga$ to get 
a compact 7-manifold $M$ with holonomy $G_2$. How can we do this?
In general we cannot, because we have no idea of how to resolve
general orbifold singularities with holonomy $G_2$. However,
suppose we can arrange that every connected component of $S$ is
locally isomorphic to either
\begin{itemize}
\setlength{\parsep}{0pt}
\setlength{\itemsep}{0pt}
\item[(a)] $T^3\t\C^2/G$, for $G$ a finite subgroup of $\SU(2)$, or
\item[(b)] ${\mathcal S}^1\t\C^3/G$, for $G$ a finite 
subgroup of $\SU(3)$ acting freely on~$\C^3\sm\{0\}$.
\end{itemize}

One can use complex algebraic geometry to find a {\it crepant
resolution} $X$ of $\C^2/G$ or $Y$ of $\C^3/G$. Then $T^3\t X$
or ${\mathcal S}^1\t Y$ gives a local model for how to resolve the
corresponding component of $S$ in $T^7/\Ga$. Thus we construct
a nonsingular, compact 7-manifold $M$ by using the patches
$T^3\t X$ or ${\mathcal S}^1\t Y$ to repair the singularities of
$T^7/\Ga$. In the case of Example \ref{ec3ex}, this means gluing
12 copies of $T^3\t X$ into $T^7/\Ga$, where $X$ is the blow-up
of $\C^2/\{\pm1\}$ at its singular point.

Now the point of using crepant resolutions is this. In both case
(a) and (b), there exists a Calabi--Yau metric on $X$ or $Y$
which is asymptotic to the flat Euclidean metric on $\C^2/G$
or $\C^3/G$. Such metrics are called {\it Asymptotically Locally
Euclidean (ALE)}. In case (a), the ALE Calabi--Yau metrics were
classified by Kronheimer \cite{Kron1,Kron2}, and exist for all
finite $G\subset\SU(2)$. In case (b), crepant resolutions of
$\C^3/G$ exist for all finite $G\subset\SU(3)$ by Roan \cite{Roan},
and the author \cite{Joyc6}, \cite[\S 8]{Joyc5} proved that they
carry ALE Calabi--Yau metrics, using a noncompact version of the
Calabi Conjecture.

By Propositions \ref{ec2prop3} and \ref{ec2prop4}, we can use
the Calabi--Yau metrics on $X$ or $Y$ to construct a torsion-free
$G_2$-structure on $T^3\t X$ or ${\mathcal S}^1\t Y$. This gives a
local model for how to resolve the singularity $T^3\t\C^2/G$ or
${\mathcal S}^1\t\C^3/G$ with holonomy $G_2$. So, this method gives
not only a way to smooth out the singularities of $T^7/\Ga$ as
a manifold, but also a family of torsion-free $G_2$-structures
on the resolution which show how to smooth out the singularities
of the $G_2$-structure.

The requirement above that $S$ be divided into connected components
of the form (a) and (b) is in fact unnecessarily restrictive. There
is a more complicated and powerful method, described in
\cite[\S 11--\S 12]{Joyc5}, for resolving singularities of a more
general kind. We require only that the singularities should
{\it locally} be of the form $\R^3\t\C^2/G$ or $\R\t\C^3/G$, for $G$
a finite subgroup of $\SU(2)$ or $\SU(3)$, and when $G\subset\SU(3)$
we do {\it not} require that $G$ act freely on~$\C^3\sm\{0\}$.

If $X$ is a crepant resolution of $\C^3/G$, where $G$ does not act
freely on $\C^3\sm\{0\}$, then the author shows \cite[\S 9]{Joyc5},
\cite{Joyc7} that $X$ carries a family of Calabi--Yau metrics
satisfying a complicated asymptotic condition at infinity, called
{\it Quasi-ALE} metrics. These yield the local models necessary to
resolve singularities locally of the form $\R\t\C^3/G$ with holonomy
$G_2$. Using this method we can resolve many orbifolds $T^7/\Ga$,
and prove the existence of large numbers of compact 7-manifolds
with holonomy~$G_2$.

\subsection{Step 3: Finding $G_2$-structures with small torsion}
\label{ec33}

For each resolution $X$ of $\C^2/G$ in case (a), and $Y$ of 
$\C^3/G$ in case (b) above, we can find a 1-parameter family 
$\{h_t:t>0\}$ of metrics with the properties
\begin{itemize}
\setlength{\parsep}{0pt}
\setlength{\itemsep}{0pt}
\item[(a)] $h_t$ is a K\"ahler metric on $X$ with $\Hol(h_t)=\SU(2)$. 
Its injectivity radius satisfies $\de(h_t)=O(t)$, its Riemann 
curvature satisfies $\bcnm{R(h_t)}{0}=O(t^{-2})$, and 
$h_t=h+O(t^4r^{-4})$ for large $r$, where $h$ is the Euclidean 
metric on $\C^2/G$, and $r$ the distance from the origin. 
\item[(b)] $h_t$ is K\"ahler on $Y$ with $\Hol(h_t)=\SU(3)$,
where $\de(h_t)=O(t)$, $\bcnm{R(h_t)}{0}=O(t^{-2})$, and
$h_t=h+O(t^6r^{-6})$ for large~$r$.
\end{itemize}
In fact we can choose $h_t$ to be isometric to $t^2h_1$, and then
(a), (b) are easy to prove.

Suppose one of the components of the singular set $S$ of $T^7/\Ga$
is locally modelled on $T^3\t\C^2/G$. Then $T^3$ has a natural flat
metric $h_{T^3}$. Let $X$ be the crepant resolution of $\C^2/G$ and
let $\{h_t:t>0\}$ satisfy property (a). Then Proposition \ref{ec2prop4}
gives a 1-parameter family of torsion-free $G_2$-structures
$(\hat\vp_t,\hat g_t)$ on $T^3\t X$ with $\hat g_t=h_{T^3}+h_t$.
Similarly, if a component of $S$ is modelled on ${\mathcal S}^1\t\C^3/G$,
using Proposition \ref{ec2prop3} we get a family of torsion-free
$G_2$-structures $(\hat\vp_t,\hat g_t)$ on~${\mathcal S}^1\t Y$.

The idea is to make a $G_2$-structure $(\vp_t,g_t)$ on $M$ by gluing 
together the torsion-free $G_2$-structures $(\hat\vp_t,\hat g_t)$ on 
the patches $T^3\t X$ and ${\mathcal S}^1\t Y$, and $(\vp_0,g_0)$ 
on $T^7/\Ga$. The gluing is done using a partition of unity. 
Naturally, the first derivative of the partition of unity introduces 
`errors', so that $(\vp_t,g_t)$ is not torsion-free. The size of the 
torsion $\nabla\vp_t$ depends on the difference $\hat\vp_t-\vp_0$ in 
the region where the partition of unity changes. On the patches 
$T^3\t X$, since $h_t-h=O(t^4r^{-4})$ and the partition of 
unity has nonzero derivative when $r=O(1)$, we find that 
$\nabla\vp_t=O(t^4)$. Similarly $\nabla\vp_t=O(t^6)$ on the patches
${\mathcal S}^1\t Y$, and so $\nabla\vp_t=O(t^4)$ on~$M$.

For small $t$, the dominant contributions to the injectivity
radius $\de(g_t)$ and Riemann curvature $R(g_t)$ are made
by those of the metrics $h_t$ on $X$ and $Y$, so we expect
$\de(g_t)=O(t)$ and $\bcnm{R(g_t)}{0}=O(t^{-2})$ by properties
(a) and (b) above. In this way we prove the following result
\cite[Th.~11.5.7]{Joyc5}, which gives the estimates on
$(\vp_t,g_t)$ that we need.

\begin{thm} On the compact\/ $7$-manifold $M$ described above,
and on many other $7$-manifolds constructed in a similar fashion,
one can write down the following data explicitly in coordinates:
\begin{itemize}
\setlength{\parsep}{0pt}
\setlength{\itemsep}{0pt}
\item Positive constants $A_1,A_2,A_3$ and\/ $\ep$,
\item A $G_2$-structure $(\vp_t,g_t)$ on $M$ with $\d\vp_t=0$
for each $t\in(0,\ep)$, and
\item A $3$-form $\psi_t$ on $M$ with\/ $\d^*\psi_t=\d^*\vp_t$
for each\/~$t\in(0,\ep)$.
\end{itemize}

These satisfy three conditions:
\begin{itemize}
\setlength{\parsep}{0pt}
\setlength{\itemsep}{0pt}
\item[{\rm(i)}] $\lnm{\psi_t}{2}\le A_1t^4$, $\cnm{\psi_t}{0}
\le A_1t^3$ and\/~$\lnm{\d^*\psi_t}{14}\le A_1t^{16/7}$,
\item[{\rm(ii)}] the injectivity radius $\de(g_t)$ satisfies
$\de(g_t)\ge A_2t$,
\item[{\rm(iii)}] the Riemann curvature $R(g_t)$ of\/ $g_t$
satisfies~$\bcnm{R(g_t)}{0}\le A_3t^{-2}$.
\end{itemize}
Here the operator $\d^*$ and the norms $\lnm{\,.\,}{2}$,
$\lnm{\,.\,}{14}$ and\/ $\cnm{\,.\,}{0}$ depend on~$g_t$.
\label{ec3thm1}
\end{thm}

Here one should regard $\psi_t$ as a {\it first integral} of 
the torsion $\nabla\vp_t$ of $(\vp_t,g_t)$. Thus the norms 
$\lnm{\psi_t}{2}$, $\cnm{\psi_t}{0}$  and $\lnm{\d^*\psi_t}{14}$
are measures of $\nabla\vp_t$. So parts (i)--(iii) say that
$\nabla\vp_t$ is small compared to the injectivity radius and
Riemann curvature of~$(M,g_t)$.

\subsection{Step 4: Deforming to a torsion-free $G_2$-structure}
\label{ec34}

We prove the following analysis result.

\begin{thm} Let\/ $A_1,A_2,A_3$ be positive constants. Then 
there exist positive constants $\ka,K$ such that whenever
$0<t\le\ka$, the following is true.

Let\/ $M$ be a compact\/ $7$-manifold, and\/ $(\vp,g)$ a 
$G_2$-structure on $M$ with\/ $\d\vp\!=\!0$. Suppose $\psi$ is 
a smooth\/ $3$-form on\/ $M$ with\/ $\d^*\psi=\d^*\vp$, and
\begin{itemize}
\setlength{\parsep}{0pt}
\setlength{\itemsep}{0pt}
\item[{\rm(i)}] $\lnm{\psi}{2}\le A_1 t^4$, 
$\cnm{\psi}{0}\le A_1 t^{1/2}$
and\/~$\lnm{\d^*\psi}{14}\le A_1$,
\item[{\rm(ii)}] the injectivity radius $\de(g)$ satisfies
$\de(g)\ge A_2 t$, and
\item[{\rm(iii)}] the Riemann curvature $R(g)$ 
satisfies~$\bcnm{R(g)}{0}\le A_3 t^{-2}$.
\end{itemize}
Then there exists a smooth, torsion-free $G_2$-structure 
$(\ti\vp,\ti g)$ on $M$ with\/~$\cnm{\ti\vp-\vp}{0}\!\le\!Kt^{1/2}$.
\label{ec3thm2}
\end{thm}

Basically, this result says that if $(\vp,g)$ is a $G_2$-structure 
on $M$, and the torsion $\nabla\vp$ is sufficiently small, then we 
can deform to a nearby $G_2$-structure $(\ti\vp,\ti g)$ that 
is torsion-free. Here is a sketch of the proof of Theorem \ref{ec3thm2}, 
ignoring several technical points. The proof is that given in 
\cite[\S 11.6--\S 11.8]{Joyc5}, which is an improved version of
the proof in~\cite{Joyc1}.

We have a 3-form $\vp$ with $\d\vp=0$ and $\d^*\vp=\d^*\psi$ for small 
$\psi$, and we wish to construct a nearby 3-form $\ti\vp$ with 
$\d\ti\vp=0$ and $\ti \d^*\ti\vp=0$. Set $\ti\vp=\vp+\d\eta$, 
where $\eta$ is a small 2-form. Then $\eta$ must satisfy a nonlinear 
p.d.e., which we write as
\e
\d^*\d\eta=-\d^*\psi+\d^*F(\d\eta),
\label{ec3eq4}
\e
where $F$ is nonlinear, satisfying~$F(\d\eta)=O\bigl(\ms{\d\eta}\bigr)$.

We solve \eq{ec3eq4} by iteration, introducing a sequence
$\{\eta_j\}_{j=0}^\infty$ with $\eta_0=0$, satisfying the
inductive equations
\e
\d^*\d\eta_{j+1}=-\d^*\psi+\d^*F(\d\eta_j),\qquad\qquad
\d^*\eta_{j+1}=0.
\label{ec3eq5}
\e
If such a sequence exists and converges to $\eta$, then taking
the limit in \eq{ec3eq5} shows that $\eta$ satisfies 
\eq{ec3eq4}, giving us the solution we want. 

The key to proving this is an {\it inductive estimate} on the
sequence $\{\eta_j\}_{j=0}^\infty$. The inductive estimate we use 
has three ingredients, the equations
\ea
\lnm{\d\eta_{j+1}}{2}&\le\lnm{\psi}{2}
+C_1\lnm{\d\eta_j}{2}\cnm{\d\eta_j}{0},
\label{ec3eq6}\\
\lnm{\nabla\d\eta_{j+1}}{14}&\le
C_2\bigl(\lnm{\d^*\psi}{14}+
\lnm{\nabla\d\eta_j}{14}\cnm{\d\eta_j}{0}
+t^{-4}\lnm{\d\eta_{j+1}}{2}\bigr),
\label{ec3eq7}\\
\cnm{\d\eta_j}{0}&\le C_3\bigl(t^{1/2}\lnm{\nabla\d\eta_j}{14}
+t^{-7/2}\lnm{\d\eta_j}{2}\bigr).
\label{ec3eq8}
\ea
Here $C_1,C_2,C_3$ are positive constants independent of $t$. Equation
\eq{ec3eq6} is obtained from \eq{ec3eq5} by taking the $L^2$-inner 
product with $\eta_{j+1}$ and integrating by parts. Using the fact that
$\d^*\vp=\d^*\psi$ and $\lnm{\psi}{2}=O(t^4)$, $\md{\psi}=O(t^{1/2})$
we get a powerful estimate of the $L^2$-norm of~$\d\eta_{j+1}$.

Equation \eq{ec3eq7} is derived from an {\it elliptic regularity 
estimate} for the operator $\d+\d^*$ acting on 3-forms on $M$. Equation
\eq{ec3eq8} follows from the {\it Sobolev embedding theorem}, 
since $L^{14}_1(M)\hookra C^0(M)$. Both \eq{ec3eq7} and
\eq{ec3eq8} are proved on small balls of radius $O(t)$ in $M$,
using parts (ii) and (iii) of Theorem \ref{ec3thm1}, and this
is where the powers of $t$ come from.

Using \eq{ec3eq6}-\eq{ec3eq8} and part (i) of Theorem \ref{ec3thm1}
we show that if 
\e
\lnm{\d\eta_j}{2}\le C_4t^4,\;\>
\lnm{\nabla\d\eta_j}{14}\le C_5,\;\>\text{and}\;\>
\cnm{\d\eta_j}{0}\le Kt^{1/2},
\label{ec3eq9}
\e
where $C_4,C_5$ and $K$ are positive constants depending on
$C_1,C_2,C_3$ and $A_1$, and if $t$ is sufficiently small, then 
the same inequalities \eq{ec3eq9} apply to $\d\eta_{j+1}$. Since 
$\eta_0=0$, by induction \eq{ec3eq9} applies for all $j$ and the 
sequence $\{\d\eta_j\}_{j=0}^\infty$ is bounded in the Banach space 
$L^{14}_1(\Lambda^3T^*M)$. One can then use standard techniques in 
analysis to prove that this sequence converges to a smooth 
limit $\d\eta$. This concludes the proof of Theorem~\ref{ec3thm2}.

\begin{figure}[htb]
{\caption{Betti numbers $(b^2,b^3)$ of compact
$G_2$-manifolds}\label{g2betti}}
{\setlength{\unitlength}{1.5pt}
\begin{picture}(215,110)(-7,-10)
\linethickness{0.25mm}
\put(0,0){\line(1,0){220}}
\put(0,0){\line(0,1){90}}
\multiput(0,15)(0,15){5}{\multiput(1,0)(1,0){220}{\vrule width .25pt
height.125pt depth.125pt}}
\multiput(10,0)(10,0){21}{\multiput(0,3)(0,3){30}{\vrule width .25pt
height.25pt depth0pt}}
\def\rla#1{\hbox to 0pt{\footnotesize\hss #1}}
\def\cla#1{\hbox to 0pt{\footnotesize\hss #1\hss}}
\put(-2,-2){\rla{0}}\put(-2,13){\rla{5}}\put(-2,28){\rla{10}}
\put(-2,43){\rla{15}}\put(-2,58){\rla{20}}\put(-2,73){\rla{25}}
\put(0,-6){\cla{0}}\put(10,-6){\cla{10}}\put(20,-6){\cla{20}}
\put(30,-6){\cla{30}}\put(40,-6){\cla{40}}\put(50,-6){\cla{50}}
\put(60,-6){\cla{60}}\put(70,-6){\cla{70}}\put(80,-6){\cla{80}}
\put(90,-6){\cla{90}}\put(100,-6){\cla{100}}\put(110,-6){\cla{110}}
\put(120,-6){\cla{120}}\put(130,-6){\cla{130}}\put(140,-6){\cla{140}}
\put(150,-6){\cla{150}}\put(160,-6){\cla{160}}\put(170,-6){\cla{170}}
\put(180,-6){\cla{180}}\put(190,-6){\cla{190}}\put(200,-6){\cla{200}}
\put(210,-6){\cla{210}}
\put(-4,95){$b^2(M)$}\put(200,6){$b^3(M)$}
\def\q#1#2{\count1=#1 \multiply\count1 by 3 \put(#2,\count1){\circle*{1.6}}}
\q{8}{47}\q{9}{46}\q{10}{45}
\q{11}{44}\q{12}{43}\q{13}{42}
\q{14}{41}\q{15}{40}\q{16}{39}
\q{4}{35}\q{5}{34}
\q{6}{33}\q{7}{32}\q{8}{31}
\q{9}{30}\q{10}{29}\q{11}{28}
\q{12}{27}\q{11}{36}
\q{5}{13}\q{3}{11}
\q{4}{17}\q{2}{10}
\q{5}{18}\q{6}{17}\q{3}{6}
\q{2}{7}\q{8}{7}\q{3}{4} 
\q{4}{99}\q{4}{95}\q{4}{91}\q{4}{87}\q{4}{83}\q{6}{89} 
\q{4}{79}\q{6}{85}\q{4}{75}\q{6}{81}\q{4}{71}\q{6}{77} 
\q{8}{83}\q{4}{67}\q{6}{73}\q{8}{79}\q{4}{63}\q{6}{69} 
\q{8}{75}\q{4}{59}\q{6}{65}\q{8}{71}\q{10}{77}\q{4}{55} 
\q{6}{61}\q{8}{67}\q{10}{73}\q{4}{51}\q{6}{57}\q{8}{63} 
\q{10}{69}\q{12}{75}\q{4}{47}\q{6}{53}\q{8}{59}\q{10}{65} 
\q{12}{71}\q{4}{43}\q{6}{49}\q{8}{55}\q{10}{61}\q{12}{67} 
\q{14}{73}\q{4}{39}\q{6}{45}\q{8}{51}\q{10}{57}\q{12}{63} 
\q{14}{69}\q{16}{75}\q{4}{35}\q{6}{41}\q{8}{47}\q{10}{53} 
\q{12}{59}\q{14}{65}\q{16}{71}\q{18}{77}\q{2}{157}\q{3}{128} 
\q{4}{99}\q{2}{153}\q{3}{124}\q{4}{95}\q{2}{149}\q{3}{120} 
\q{4}{91}\q{2}{145}\q{3}{116}\q{4}{87}\q{2}{141}\q{4}{147} 
\q{3}{112}\q{5}{118}\q{4}{83}\q{6}{89}\q{2}{137}\q{4}{143} 
\q{3}{108}\q{5}{114}\q{4}{79}\q{6}{85}\q{2}{133}\q{4}{139} 
\q{3}{104}\q{5}{110}\q{4}{75}\q{6}{81}\q{2}{129}\q{4}{135} 
\q{6}{141}\q{3}{100}\q{5}{106}\q{7}{112}\q{4}{71}\q{6}{77} 
\q{8}{83}\q{2}{125}\q{4}{131}\q{6}{137}\q{3}{96}\q{5}{102} 
\q{7}{108}\q{4}{67}\q{6}{73}\q{8}{79}\q{2}{121}\q{4}{127} 
\q{6}{133}\q{3}{92}\q{5}{98}\q{7}{104}\q{4}{63}\q{6}{69} 
\q{8}{75}\q{4}{123}\q{6}{129}\q{8}{135}\q{5}{94}\q{7}{100} 
\q{9}{106}\q{6}{65}\q{8}{71}\q{10}{77}\q{6}{125}\q{8}{131} 
\q{7}{96}\q{9}{102}\q{8}{67}\q{10}{73}\q{8}{127}\q{10}{133} 
\q{9}{98}\q{11}{104}\q{10}{69}\q{12}{75}\q{10}{129}\q{11}{100} 
\q{12}{71}\q{12}{131}\q{13}{102}\q{14}{73}\q{14}{133}\q{15}{104} 
\q{16}{75}\q{16}{135}\q{17}{106}\q{18}{77}\q{0}{215}\q{1}{186} 
\q{2}{157}\q{3}{128}\q{4}{99}
\q{5}{74}\q{5}{70}\q{5}{66}\q{5}{62}\q{5}{58}\q{7}{64}
\q{6}{65}\q{6}{61}\q{6}{57}\q{6}{53}\q{6}{49}\q{8}{55}
\q{7}{56}\q{7}{52}\q{7}{48}\q{7}{44}\q{7}{40}\q{9}{46}
\q{9}{46}\q{8}{47}\q{8}{43}\q{8}{39}\q{8}{35}\q{8}{31}
\q{10}{37}\q{10}{53}\q{10}{49}\q{10}{45}\q{10}{41}
\q{10}{37}\q{12}{43}\q{12}{55}\q{12}{51}\q{12}{47}
\q{12}{43}\q{11}{44}\q{11}{40}\q{11}{36}\q{11}{32}
\q{11}{28}\q{13}{34}\q{13}{46}\q{13}{42}\q{13}{38}
\q{13}{34}\q{15}{40}\q{14}{41}\q{14}{37}\q{14}{33}
\q{14}{29}\q{14}{25}\q{16}{31}\q{16}{43}\q{16}{39}
\q{16}{35}\q{16}{31}\q{18}{45}\q{18}{41}\q{18}{37}
\q{20}{43}\q{17}{38}\q{17}{34}\q{17}{30}\q{17}{26}
\q{19}{28}\q{19}{40}\q{19}{36}\q{19}{32}\q{19}{28}
\q{21}{42}\q{21}{38}\q{21}{34}\q{23}{44}\q{23}{40}
\q{25}{46}\q{3}{84}\q{3}{80}\q{3}{76}\q{3}{72}\q{3}{68}
\q{4}{79}\q{4}{75}\q{4}{71}\q{4}{67}\q{4}{63}\q{4}{75}
\q{4}{71}\q{4}{67}\q{4}{63}\q{4}{59}\q{5}{70}\q{5}{66}
\q{5}{62}\q{5}{58}\q{5}{54}\q{5}{66}\q{5}{62}\q{5}{58}
\q{5}{54}\q{5}{50}\q{6}{61}\q{6}{57}\q{6}{53}\q{6}{49}
\q{6}{45}\q{7}{56}\q{8}{51}\q{6}{57}\q{6}{53}\q{6}{49}
\q{6}{45}\q{7}{52}\q{7}{48}\q{7}{44}\q{7}{40}\q{8}{63}
\q{8}{59}\q{8}{55}\q{8}{51}\q{8}{47}\q{9}{58}\q{9}{54}
\q{9}{50}\q{9}{46}\q{9}{42}\q{10}{65}\q{10}{61}\q{10}{57}
\q{10}{53}\q{11}{60}\q{11}{56}\q{11}{52}\q{11}{48}\q{9}{54}
\q{9}{50}\q{9}{46}\q{10}{49}\q{10}{45}\q{10}{41}\q{11}{56}
\q{11}{52}\q{11}{48}\q{12}{51}\q{12}{47}\q{12}{43}\q{13}{50}
\q{14}{45}\q{12}{51}\q{12}{47}\q{13}{46}\q{13}{42}\q{14}{53}
\q{14}{49}\q{15}{48}\q{15}{44}\q{16}{55}\q{16}{51}\q{17}{50}
\q{17}{46}\q{18}{53}\q{19}{48}\q{15}{48}\q{16}{43}\q{17}{50}
\q{18}{45}\q{19}{52}\q{20}{47}\q{21}{54}\q{22}{49}\q{23}{56}
\q{24}{51}\q{1}{142}\q{2}{137}\q{3}{132}\q{3}{128}\q{3}{124}
\q{3}{120}\q{5}{122}\q{2}{113}\q{3}{108}\q{4}{103}\q{4}{99}
\q{4}{95}\q{4}{91}\q{6}{93}\q{3}{84}\q{4}{79}\q{5}{74}
\q{5}{70}\q{5}{66}\q{5}{62}\q{7}{64}
\q{28}{43}
\q{5}{26}\q{10}{29}\q{15}{32}\q{20}{35}
\q{10}{13}\q{5}{10}
\q{2}{19}\q{3}{19}\q{4}{19}\q{5}{19}\q{6}{19}\q{7}{19}\q{8}{19}
\q{9}{19}\q{10}{19}\q{11}{19}
\end{picture}}
\end{figure}

From Theorems \ref{ec3thm1} and \ref{ec3thm2} we see that the compact
7-manifold $M$ constructed in Step 2 admits torsion-free
$G_2$-structures $(\ti\vp,\ti g)$. Theorem \ref{ec2thm3} then
shows that $\Hol(\ti g)=G_2$ if and only if $\pi_1(M)$ is finite.
In the example above $M$ is simply-connected, and so $\pi_1(M)=\{1\}$
and $M$ has metrics with holonomy $G_2$, as we want.

By considering different groups $\Ga$ acting on $T^7$, and also 
by finding topologically distinct resolutions $M_1,\dots,M_k$ of 
the same orbifold $T^7/\Ga$, we can construct many compact 
Riemannian 7-manifolds with holonomy $G_2$. A good number of
examples are given in \cite[\S 12]{Joyc5}. Figure \ref{g2betti}
displays the Betti numbers of compact, simply-connected 7-manifolds
with holonomy $G_2$ constructed there. There are 252 different sets
of Betti numbers.

Examples are also known \cite[\S 12.4]{Joyc5} of compact 7-manifolds
with holonomy $G_2$ with finite, nontrivial fundamental group. It
seems likely to the author that the Betti numbers given in Figure
\ref{g2betti} are only a small proportion of the Betti numbers of
all compact, simply-connected 7-manifolds with holonomy~$G_2$.

\subsection{Other constructions of compact $G_2$-manifolds}
\label{ec35}

Here are two other methods, taken from \cite[\S 11.9]{Joyc5},
of constructing compact 7-manifolds with holonomy $G_2$. The
first was outlined by the author in~\cite[\S 4.3]{Joyc2}.
\medskip

\noindent{\bf Method 1.} Let $(Y,J,h)$ be a Calabi--Yau 3-fold, 
with K\"ahler form $\om$ and holomorphic volume form $\th$. Suppose
$\si:Y\ra Y$ is an involution, satisfying $\si^*(h)=h$,
$\si^*(J)=-J$ and $\si^*(\th)=\bar\th$. We call
$\si$ a {\it real structure} on $Y$. Let $N$ be the fixed point 
set of $\si$ in $Y$. Then $N$ is a real 3-dimensional submanifold 
of $Y$, and is in fact a special Lagrangian 3-fold.

Let ${\mathcal S}^1=\R/\Z$, and define a torsion-free $G_2$-structure 
$(\vp,g)$ on ${\mathcal S}^1\t Y$ as in Proposition \ref{ec2prop3}. 
Then $\vp=\d x\w\om+\Re\th$, where $x\in\R/\Z$ is the coordinate 
on ${\mathcal S}^1$. Define $\hat\si:{\mathcal S}^1\t Y\ra{\mathcal S}^1\t Y$ 
by $\hat\si\bigl((x,y)\bigr)=\bigl(-x,\si(y)\bigr)$. Then $\hat\si$
preserves $(\vp,g)$ and $\hat\si^2=1$. The fixed points of 
$\hat\si$ in ${\mathcal S}^1\t Y$ are $\{\Z,\ha+\Z\}\t N$. Thus 
$({\mathcal S}^1\t Y)/\an{\hat\si}$ is an orbifold. Its singular
set is 2 copies of $N$, and each singular point is modelled 
on~$\R^3\t\R^4/\{\pm 1\}$.

We aim to resolve $({\mathcal S}^1\t Y)/\an{\hat\si}$ to get a
compact 7-manifold $M$ with holonomy $G_2$. Locally, each
singular point should be resolved like $\R^3\t X$, where
$X$ is an ALE Calabi--Yau 2-fold asymptotic to $\C^2/\{\pm 1\}$.
There is a 3-dimensional family of such $X$, and we need to
choose one member of this family for each singular point in
the singular set.

Calculations by the author indicate that the data needed to do 
this is a closed, coclosed 1-form $\al$ on $N$ that is nonzero
at every point of $N$. The existence of a suitable 1-form $\al$
depends on the metric on $N$, which is the restriction of the
metric $g$ on $Y$. But $g$ comes from the solution of the Calabi
Conjecture, so we know little about it. This may make the method
difficult to apply in practice.
\medskip

The second method has been successfully applied by Kovalev
\cite{Kova}, and is based on an idea due to Simon Donaldson.
\medskip

\noindent{\bf Method 2.} Let $X$ be a projective complex 3-fold
with canonical bundle $K_X$, and $s$ a holomorphic section of 
$K_X^{-1}$ which vanishes to order 1 on a smooth divisor $D$ 
in $X$. Then $D$ has trivial canonical bundle, so $D$ is $T^4$
or $K3$. Suppose $D$ is a $K3$ surface. Define $Y=X\sm D$,
and suppose $Y$ is simply-connected.

Then $Y$ is a noncompact complex 3-fold with $K_Y$ trivial, 
and one infinite end modelled on $D\t{\mathcal S}^1\t[0,\iy)$.
Using a version of the proof of the Calabi Conjecture for
noncompact manifolds one constructs a complete Calabi--Yau
metric $h$ on $Y$, which is asymptotic to the product on
$D\t{\mathcal S}^1\t[0,\iy)$ of a Calabi--Yau metric on $D$,
and Euclidean metrics on ${\mathcal S}^1$ and $[0,\iy)$. We
call such metrics {\it Asymptotically Cylindrical}.

Suppose we have such a metric on $Y$. Define a torsion-free 
$G_2$-structure $(\vp,g)$ on ${\mathcal S}^1\t Y$ as in Proposition
\ref{ec2prop3}. Then ${\mathcal S}^1\t Y$ is a noncompact $G_2$-manifold
with one end modelled on $D\t T^2\t[0,\iy)$, whose metric is
asymptotic to the product on $D\t T^2\t[0,\iy)$ of a Calabi--Yau
metric on $D$, and Euclidean metrics on $T^2$ and~$[0,\iy)$. 

Donaldson and Kovalev's idea is to take two such products
${\mathcal S}^1\t Y_1$ and ${\mathcal S}^1\t Y_2$ whose infinite ends 
are isomorphic in a suitable way, and glue them together to get 
a compact 7-manifold $M$ with holonomy $G_2$. The gluing process
swaps round the ${\mathcal S}^1$ factors. That is, the ${\mathcal S}^1$
factor in ${\mathcal S}^1\t Y_1$ is identified with the asymptotic
${\mathcal S}^1$ factor in $Y_2\sim D_2\t{\mathcal S}^1\t[0,\iy)$, and
vice versa.

\section{Compact $\Spin(7)$-manifolds from Calabi--Yau 4-orbifolds}
\label{ec4}

In a very similar way to the $G_2$ case, one can construct examples
of compact 8-manifolds with holonomy $\Spin(7)$ by resolving the
singularities of torus orbifolds $T^8/\Ga$. This is done in
\cite{Joyc3} and \cite[\S 13--\S 14]{Joyc5}. In \cite[\S 14]{Joyc5},
examples are constructed which realize 181 different sets of Betti
numbers. Two compact 8-manifolds with holonomy $\Spin(7)$ and the
same Betti numbers may be distinguished by the cup products on their
cohomologies (examples of this are given in \cite[\S 3.4]{Joyc3}),
so they probably represent rather more than 181 topologically
distinct 8-manifolds.

The main differences with the $G_2$ case are, firstly, that the
technical details of the analysis are different and harder, and
secondly, that the singularities that arise are typically more
complicated and more tricky to resolve. One reason for this is that
in the $G_2$ case the singular set is made up of 1 and 3-dimensional
pieces in a 7-dimensional space, so one can often arrange for the
pieces to avoid each other, and resolve them independently.

But in the $\Spin(7)$ case the singular set is typically made up
of 4-dimensional pieces in an 8-dimensional space, so they nearly
always intersect. There are also topological constraints arising
from the $\hat A$-genus, which do not apply in the $G_2$ case.
The moral appears to be that when you increase the dimension,
things become more difficult.

Anyway, we will not discuss this further, as the principles are
very similar to the $G_2$ case above. Instead, we will discuss
an entirely different construction of compact 8-manifolds with
holonomy $\Spin(7)$ developed by the author in \cite{Joyc4} and
\cite[\S 15]{Joyc5}, a little like Method 1 of \S\ref{ec35}. In
this we start from a {\it Calabi--Yau $4$-orbifold} rather than
from $T^8$. The construction can be divided into five steps.
\begin{list}{}{\setlength{\leftmargin}{40pt}
\setlength{\labelwidth}{40pt}}
\item[Step 1.] Find a compact, complex 4-orbifold $(Y,J)$
satisfying the conditions:
\begin{itemize}
\setlength{\parsep}{0pt}
\setlength{\itemsep}{0pt}
\item[(a)] $Y$ has only finitely many singular points
$p_1,\ldots,p_k$, for~$k\ge 1$.
\item[(b)] $Y$ is modelled on $\C^4/\an{i}$ near each $p_j$,
where $i$ acts on $\C^4$ by complex multiplication.
\item[(c)] There exists an antiholomorphic involution
$\si:Y\ra Y$ whose fixed point set is~$\{p_1,\ldots,p_k\}$.
\item[(d)] $Y\sm\{p_1,\ldots,p_k\}$ is simply-connected,
and~$h^{2,0}(Y)=0$.
\end{itemize}
\item[Step 2.] Choose a $\si$-invariant K\"ahler class on $Y$.
Then by Theorem \ref{ec2thm5} there exists a unique $\si$-invariant
Ricci-flat K\"ahler metric $g$ in this K\"ahler class. Let $\om$
be the K\"ahler form of $g$. Let $\th$ be a holomorphic volume
form for $(Y,J,g)$. By multiplying $\th$ by ${\rm e}^{i\phi}$ if
necessary, we can arrange that~$\si^*(\th)=\bar\th$. 

Define $\Om=\ha\om\w\om+\Re\th$. Then $(\Om,g)$ is a torsion-free
$\Spin(7)$-structure on $Y$, by Proposition \ref{ec2prop5}. Also,
$(\Om,g)$ is $\si$-invariant, as $\si^*(\om)=-\om$ and $\si^*(\th)
=\bar\th$. Define $Z=Y/\an{\si}$. Then $Z$ is a compact real
8-orbifold with isolated singular points $p_1,\ldots,p_k$, and
$(\Om,g)$ pushes down to a torsion-free $\Spin(7)$-structure
$(\Om,g)$ on $Z$.
\item[Step 3.] $Z$ is modelled on $\R^8/G$ near each $p_j$, where
$G$ is a certain finite subgroup of $\Spin(7)$ with $\md{G}=8$.
We can write down two explicit, topologically distinct ALE
$\Spin(7)$-manifolds $X_1,X_2$ asymptotic to $\R^8/G$. Each
carries a 1-parameter family of homothetic ALE metrics $h_t$
for $t>0$ with $\Hol(h_t)=\Z_2\lt\SU(4)\subset\Spin(7)$.

For $j=1,\ldots,k$ we choose $i_j=1$ or 2, and resolve the
singularities of $Z$ by gluing in $X_{i_j}$ at the singular
point $p_j$ for $j=1,\ldots,k$, to get a compact, nonsingular
8-manifold $M$, with projection~$\pi:M\ra Z$.

\item[Step 4.] On $M$, we explicitly write down a 1-parameter 
family of $\Spin(7)$-structures $(\Om_t,g_t)$ depending on $t\in
(0,\ep)$. They are not torsion-free, but have small torsion when
$t$ is small. As $t\ra 0$, the $\Spin(7)$-structure $(\Om_t,g_t)$
converges to the singular $\Spin(7)$-structure~$\pi^*(\Om_0,g_0)$.
\item[Step 5.] We prove using analysis that for sufficiently
small $t$, the $\Spin(7)$-structure $(\Om_t,g_t)$ on $M$, with
small torsion, can be deformed to a $\Spin(7)$-structure
$(\ti\Om_t,\ti g_t)$, with zero torsion.

It turns out that if $i_j=1$ for $j=1,\ldots,k$ we have
$\pi_1(M)\cong\Z_2$ and $\Hol(\ti g_t)=\Z_2\lt\SU(4)$, and for
the other $2^k-1$ choices of $i_1,\ldots,i_k$ we have $\pi_1(M)=\{1\}$
and $\Hol(\ti g_t)=\Spin(7)$. So $\ti g_t$ is a metric with holonomy
$\Spin(7)$ on the compact 8-manifold $M$ for~$(i_1,\ldots,i_k)\ne
(1,\ldots,1)$.
\end{list}

Once we have completed Step 1, Step 2 is immediate. Steps 4 and 5
are analogous to Steps 3 and 4 of \S\ref{ec3}, and can be done using
the techniques and analytic results developed by the author for
the first $T^8/\Ga$ construction of compact $\Spin(7)$-manifolds,
\cite{Joyc3}, \cite[\S 13]{Joyc5}. So the really new material is
in Steps 1 and 3, and we will discuss only these.

\subsection{Step 1: An example}
\label{ec41}

We do Step 1 using complex algebraic geometry. The problem is that
conditions (a)--(d) above are very restrictive, so it is not that
easy to find {\it any} $Y$ satisfying all four conditions. All the
examples $Y$ the author has found are constructed using {\it weighted
projective spaces}, an important class of complex orbifolds.

\begin{dfn} Let $m\ge 1$ be an integer, and $a_0,a_1,\ldots,a_m$
positive integers with highest common factor 1. Let $\C^{m+1}$ have 
complex coordinates on $(z_0,\ldots,z_m)$, and define an action of 
the complex Lie group $\C^*$ on $\C^{m+1}$ by
\begin{equation*}
(z_0,\ldots,z_m)\,{\buildrel u\over\longmapsto}
(u^{a_0}z_0,\ldots,u^{a_m}z_m),\qquad\text{for $u\in\C^*$.}
\end{equation*}
The {\it weighted projective space} $\CP^m_{a_0,\ldots,a_m}$
is $\bigl(\C^{m+1}\sm\{0\}\bigr)/\C^*$. The $\C^*$-orbit of
$(z_0,\ldots,z_m)$ is written~$[z_0,\ldots,z_m]$.
\label{ec4def}
\end{dfn}

Here is the simplest example the author knows.

\begin{ex} Let $Y$ be the hypersurface of degree 12 in 
$\CP^5_{1,1,1,1,4,4}$ given by
\begin{equation*}
Y=\bigl\{[z_0,\ldots,z_5]\in\CP^5_{1,1,1,1,4,4}:z_0^{12}+
z_1^{12}+z_2^{12}+z_3^{12}+z_4^3+z_5^3=0\bigr\}.
\end{equation*}
Calculation shows that $Y$ has trivial canonical bundle and three 
singular points $p_1\!=\![0,0,0,0,1,-1]$, $p_2\!=\![0,0,0,0,1,e^{\pi i/3}]$ 
and $p_3\!=\![0,0,0,0,1,e^{-\pi i/3}]$, modelled on~$\C^4/\an{i}$.

Now define a map $\si:Y\ra Y$ by
\begin{equation*}
\si:[z_0,\ldots,z_5]\longmapsto[\bar z_1,-\bar z_0,\bar z_3,
-\bar z_2,\bar z_5,\bar z_4].
\end{equation*}
Note that $\si^2=1$, though this is not immediately obvious,
because of the geometry of $\CP^5_{1,1,1,1,4,4}$. It can be shown 
that conditions (a)--(d) of Step 1 above hold for $Y$ and~$\si$.
\label{ec4ex1}
\end{ex}

More suitable 4-folds $Y$ may be found by taking hypersurfaces
or complete intersections in other weighted projective spaces,
possibly also dividing by a finite group, and then doing a
crepant resolution to get rid of any singularities that we don't
want. Examples are given in \cite{Joyc4}, \cite[\S 15]{Joyc5}.

\subsection{Step 3: Resolving $\R^8/G$}
\label{ec42}

Define $\al,\be:\R^8\ra\R^8$ by
\begin{equation*}
\begin{split}
\al:(x_1,\ldots,x_8)&\mapsto(-x_2,x_1,-x_4,x_3,-x_6,x_5,-x_8,x_7),\\
\be:(x_1,\ldots,x_8)&\mapsto(x_3,-x_4,-x_1,x_2,x_7,-x_8,-x_5,x_6).
\end{split}
\end{equation*}
Then $\al,\be$ preserve $\Om_0$ given in \eq{ec2eq2}, so they lie in
$\Spin(7)$. Also $\al^4=\be^4=1$, $\al^2=\be^2$ and $\al\be=\be\al^3$.
Let $G=\an{\al,\be}$. Then $G$ is a finite nonabelian subgroup of
$\Spin(7)$ of order 8, which acts freely on $\R^8\sm\{0\}$. One
can show that if $Z$ is the compact $\Spin(7)$-orbifold constructed
in Step 2 above, then $T_{p_j}Z$ is isomorphic to $\R^8/G$ for
$j=1,\ldots,k$, with an isomorphism identifying the
$\Spin(7)$-structures $(\Om,g)$ on $Z$ and $(\Om_0,g_0)$ on $\R^8/G$,
such that $\be$ corresponds to the $\si$-action on~$Y$.

In the next two examples we shall construct two different ALE
$\Spin(7)$-manifolds $(X_1,\Om_1,g_1)$ and $(X_2,\Om_2,g_2)$
asymptotic to~$\R^8/G$. 

\begin{ex} Define complex coordinates $(z_1,\ldots,z_4)$ on $\R^8$ by
\begin{equation*}
(z_1,z_2,z_3,z_4)=(x_1+ix_2,x_3+ix_4,x_5+ix_6,x_7+ix_8),
\end{equation*}
Then $g_0=\ms{\d z_1}+\cdots+\ms{\d z_4}$, and $\Om_0=\ha\om_0
\w\om_0+\Re(\th_0)$, where $\om_0$ and $\th_0$ are the 
usual K\"ahler form and complex volume form on $\C^4$. In these 
coordinates, $\al$ and $\be$ are given by
\begin{equation}
\begin{split}
\al:(z_1,\ldots,z_4)&\mapsto(iz_1,iz_2,iz_3,iz_4),\\
\be:(z_1,\ldots,z_4)&\mapsto(\bar z_2,-\bar z_1,\bar z_4,-\bar z_3).
\end{split}
\label{ec4eq1}
\end{equation}

Now $\C^4/\an{\al}$ is a complex singularity, as $\al\in\SU(4)$.
Let $(Y_1,\pi_1)$ be the blow-up of $\C^4/\an{\al}$ at 0. Then $Y_1$
is the unique crepant resolution of $\C^4/\an{\al}$. The action of
$\be$ on $\C^4/\an{\al}$ lifts to a {\it free} antiholomorphic
map $\be:Y_1\ra Y_1$ with $\be^2=1$. Define $X_1=Y_1/\an{\be}$. 
Then $X_1$ is a nonsingular 8-manifold, and the projection 
$\pi_1:Y_1\ra\C^4/\an{\al}$ pushes down to~$\pi_1:X_1\ra\R^8/G$. 

There exist ALE Calabi--Yau metrics $g_1$ on $Y_1$, which were
written down explicitly by Calabi \cite[p.~285]{Cal}, and are
invariant under the action of $\be$ on $Y_1$. Let $\om_1$ be
the K\"ahler form of $g_1$, and $\th_1=\pi_1^*(\th_0)$ the
holomorphic volume form on $Y_1$. Define $\Om_1=\ha\om_1\w
\om_1+\Re(\th_1)$. Then $(\Om_1,g_1)$ is a torsion-free 
$\Spin(7)$-structure on $Y_1$, as in Proposition~\ref{ec2prop5}.

As $\be^*(\om_1)=-\om_1$ and $\be^*(\th_1)=\bar\th_1$, we 
see that $\be$ preserves $(\Om_1,g_1)$. Thus $(\Om_1,g_1)$ pushes 
down to a torsion-free $\Spin(7)$-structure $(\Om_1,g_1)$ on 
$X_1$. Then $(X_1,\Om_1,g_1)$ is an {\it ALE\/ $\Spin(7)$-manifold} 
asymptotic to~$\R^8/G$.
\label{ec4ex2}
\end{ex}

\begin{ex} Define new complex coordinates $(w_1,\ldots,w_4)$ on $\R^8$ by
\begin{equation*}
(w_1,w_2,w_3,w_4)=(-x_1+ix_3,x_2+ix_4,-x_5+ix_7,x_6+ix_8).
\end{equation*}
Again we find that $g_0=\ms{\d w_1}+\cdots+\ms{\d w_4}$ and 
$\Om_0=\ha\om_0\w\om_0+\Re(\th_0)$. In these 
coordinates, $\al$ and $\be$ are given by
\begin{equation}
\begin{split}
\al:(w_1,\ldots,w_4)&\mapsto(\bar w_2,-\bar w_1,\bar w_4,-\bar w_3),\\
\be:(w_1,\ldots,w_4)&\mapsto(iw_1,iw_2,iw_3,iw_4).
\end{split}
\label{ec4eq2}
\end{equation}
Observe that \eq{ec4eq1} and \eq{ec4eq2} are the same, except that 
the r\^oles of $\al,\be$ are reversed. Therefore we can use 
the ideas of Example \ref{ec4ex2} again.

Let $Y_2$ be the crepant resolution of $\C^4/\an{\be}$. The action 
of $\al$ on $\C^4/\an{\be}$ lifts to a free antiholomorphic
involution of $Y_2$. Let $X_2=Y_2/\an{\al}$. Then $X_2$ is nonsingular, 
and carries a torsion-free $\Spin(7)$-structure $(\Om_2,g_2)$, making 
$(X_2,\Om_2,g_2)$ into an ALE $\Spin(7)$-manifold asymptotic to~$\R^8/G$.
\label{ec4ex3}
\end{ex}

We can now explain the remarks on holonomy groups at the end of
Step 5. The holonomy groups $\Hol(g_i)$ of the metrics $g_1,g_2$
in Examples \ref{ec4ex2} and \ref{ec4ex3} are both isomorphic to
$\Z_2\lt\SU(4)$, a subgroup of $\Spin(7)$. However, they are two
{\it different} inclusions of $\Z_2\lt\SU(4)$ in $\Spin(7)$, as in
the first case the complex structure is $\al$ and in the second~$\be$.

The $\Spin(7)$-structure $(\Om,g)$ on $Z$ also has holonomy
$\Hol(g)=\Z_2\lt\SU(4)$. Under the natural identifications we have
$\Hol(g_1)=\Hol(g)$ but $\Hol(g_2)\ne\Hol(g)$ as subgroups of
$\Spin(7)$. Therefore, if we choose $i_j=1$ for all $j=1,\ldots,k$,
then $Z$ and $X_{i_j}$ all have the same holonomy group
$\Z_2\lt\SU(4)$, so they combine to give metrics $\ti g_t$ on
$M$ with~$\Hol(\ti g_t)=\Z_2\lt\SU(4)$.

However, if $i_j=2$ for some $j$ then the holonomy of $g$ on $Z$
and $g_{i_j}$ on $X_{i_j}$ are {\it different} $\Z_2\lt\SU(4)$
subgroups of $\Spin(7)$, which together generate the whole group
$\Spin(7)$. Thus they combine to give metrics $\ti g_t$ on $M$
with~$\Hol(\ti g_t)=\Spin(7)$.

\subsection{Conclusions}
\label{ec43}

The author was able in \cite{Joyc4} and \cite[Ch.~15]{Joyc5} to
construct compact 8-manifolds with holonomy $\Spin(7)$ realizing
14 distinct sets of Betti numbers, which are given in Table
\ref{ec7betti}. Probably there are many other examples which can
be produced by similar methods.

\begin{table}[htb]
\centering
{\caption{Betti numbers $(b^2,b^3,b^4)$ of compact
$\Spin(7)$-manifolds}\label{ec7betti}}
{\begin{tabular}{ccccc}
\hline
\vphantom{$\bigr)^{k^k}$}
(4,\,33,\,200) & (3,\,33,\,202) & (2,\,33,\,204) & 
(1,\,33,\,206) & (0,\,33,\,208) \\ 
(1,\,0,\,908)  & (0,\,0,\,910)  & (1,\,0,\,1292) & 
(0,\,0,\,1294) & (1,\,0,\,2444) \\ 
(0,\,0,\,2446) & (0,\,6,\,3730) & (0,\,0,\,4750) & 
(0,\,0,\,11\,662)\vphantom{$\bigr)_{p_p}$} \\
\hline
\end{tabular}}
\end{table}

Comparing these Betti numbers with those of the compact 8-manifolds
constructed in \cite[Ch.~14]{Joyc5} by resolving torus orbifolds
$T^8/\Gamma$, we see that these examples the middle Betti number
$b^4$ is much bigger, as much as $11\,662$ in one case.

Given that the two constructions of compact 8-manifolds with holonomy 
$\Spin(7)$ that we know appear to produce sets of 8-manifolds with 
rather different `geography', it is tempting to speculate that the 
set of all compact 8-manifolds with holonomy $\Spin(7)$ may be rather
large, and that those constructed so far are a small sample with 
atypical behaviour.
\vfil\eject

\centerline{\LARGE\scshape Part II. Calibrated Geometry}
\medskip

\section{Introduction to calibrated geometry}
\label{ec5}

{\it Calibrated geometry} was introduced in the seminal paper
of Harvey and Lawson \cite{HaLa}. We introduce the basic
ideas in \S\ref{ec51}--\S\ref{ec52}, and then discuss the
$G_2$ calibrations in more detail in \S\ref{ec53}--\S\ref{ec55},
and the $\Spin(7)$ calibration in~\S\ref{ec56}.

\subsection{Calibrations and calibrated submanifolds}
\label{ec51}

We begin by defining {\it calibrations} and {\it calibrated 
submanifolds}, following Harvey and Lawson~\cite{HaLa}.

\begin{dfn} Let $(M,g)$ be a Riemannian manifold. An {\it oriented
tangent\/ $k$-plane} $V$ on $M$ is a vector subspace $V$ of
some tangent space $T_xM$ to $M$ with $\dim V=k$, equipped
with an orientation. If $V$ is an oriented tangent $k$-plane
on $M$ then $g\vert_V$ is a Euclidean metric on $V$, so 
combining $g\vert_V$ with the orientation on $V$ gives a 
natural {\it volume form} $\vol_V$ on $V$, which is a 
$k$-form on~$V$.

Now let $\vp$ be a closed $k$-form on $M$. We say that
$\vp$ is a {\it calibration} on $M$ if for every oriented
$k$-plane $V$ on $M$ we have $\vp\vert_V\le \vol_V$. Here
$\vp\vert_V=\al\cdot\vol_V$ for some $\al\in\R$, and 
$\vp\vert_V\le\vol_V$ if $\al\le 1$. Let $N$ be an 
oriented submanifold of $M$ with dimension $k$. Then 
each tangent space $T_xN$ for $x\in N$ is an oriented
tangent $k$-plane. We say that $N$ is a {\it calibrated 
submanifold\/} if $\vp\vert_{T_xN}=\vol_{T_xN}$ for all~$x\in N$.
\label{ec5def1}
\end{dfn}

It is easy to show that calibrated submanifolds are automatically
{\it minimal submanifolds} \cite[Th.~II.4.2]{HaLa}. We prove this
in the compact case, but noncompact calibrated submanifolds are
locally volume-minimizing as well.

\begin{prop} Let\/ $(M,g)$ be a Riemannian manifold, $\vp$ a 
calibration on $M$, and\/ $N$ a compact $\vp$-submanifold 
in $M$. Then $N$ is volume-minimizing in its homology class.
\label{ec5prop1}
\end{prop}

\begin{proof} Let $\dim N=k$, and let $[N]\in H_k(M,\R)$ 
and $[\vp]\in H^k(M,\R)$ be the homology and cohomology 
classes of $N$ and $\vp$. Then
\begin{equation*}
[\vp]\cdot[N]=\int_{x\in N}\vp\big\vert_{T_xN}=
\int_{x\in N}{\ts\vol_{T_xN}}=\Vol(N),
\end{equation*}
since $\vp\vert_{T_xN}=\vol_{T_xN}$ for each $x\in N$, as
$N$ is a calibrated submanifold. If $N'$ is any other compact 
$k$-submanifold of $M$ with $[N']=[N]$ in $H_k(M,\R)$, then
\begin{equation*}
[\vp]\cdot[N]=[\vp]\cdot[N']=\int_{x\in N'}\vp\big\vert_{T_xN'}
\le\int_{x\in N'}{\ts\vol_{T_xN'}}=\Vol(N'),
\end{equation*}
since $\vp\vert_{T_xN'}\le\vol_{T_xN'}$ because $\vp$ is a 
calibration. The last two equations give $\Vol(N)\le\Vol(N')$. 
Thus $N$ is volume-minimizing in its homology class.
\end{proof}

Now let $(M,g)$ be a Riemannian manifold with a calibration $\vp$,
and let $\iota:N\ra M$ be an immersed submanifold. Whether
$N$ is a $\vp$-submanifold depends upon the tangent spaces of 
$N$. That is, it depends on $\iota$ and its first derivative. 
So, to be calibrated with respect to $\vp$ is a {\it first-order}
partial differential equation on $\iota$. But if $N$ is calibrated
then $N$ is minimal, and to be minimal is a {\it second-order}
partial differential equation on~$\iota$.

One moral is that the calibrated equations, being first-order,
are often easier to solve than the minimal submanifold equations,
which are second-order. So calibrated geometry is a fertile source
of examples of minimal submanifolds. 

\subsection{Calibrated submanifolds and special holonomy}
\label{ec52}

Next we explain the connection with Riemannian holonomy. Let 
$G\subset{\rm O}(n)$ be a possible holonomy group of a Riemannian 
metric. In particular, we can take $G$ to be one of the
holonomy groups $\U(m)$, $\SU(m)$, $\Sp(m)$, $G_2$ or 
Spin(7) from Berger's classification. Then $G$ acts on 
the $k$-forms $\La^k(\R^n)^*$ on $\R^n$, so we can look
for $G$-invariant $k$-forms on~$\R^n$.

Suppose $\vp_0$ is a nonzero, $G$-invariant $k$-form on 
$\R^n$. By rescaling $\vp_0$ we can arrange that for each 
oriented $k$-plane $U\subset\R^n$ we have $\vp_0\vert_U\le\vol_U$, 
and that $\vp_0\vert_U=\vol_U$ for at least one such $U$. Then
$\vp_0\vert_{\ga\cdot U}=\vol_{\ga\cdot U}$ by $G$-invariance,
so $\ga\cdot U$ is a calibrated $k$-plane for all $\ga\in G$.
Thus the family of $\vp_0$-calibrated $k$-planes in $\R^n$
is reasonably large, and it is likely the calibrated
submanifolds will have an interesting geometry.

Now let $M$ be a manifold of dimension $n$, and $g$ a metric
on $M$ with Levi-Civita connection $\nabla$ and holonomy 
group $G$. Then by Theorem \ref{ec2thm1} there is a $k$-form 
$\vp$ on $M$ with $\nabla\vp=0$, corresponding to $\vp_0$. 
Hence $\d\vp=0$, and $\vp$ is closed. Also, the condition 
$\vp_0\vert_U\le\vol_U$ for all oriented $k$-planes $U$ in $\R^n$ 
implies that $\vp\vert_V\le\vol_V$ for all oriented tangent 
$k$-planes $V$ in $M$. Thus $\vp$ is a {\it calibration} on~$M$. 

This gives us a general method for finding interesting 
calibrations on manifolds with reduced holonomy. Here are 
the most significant examples of this.
\begin{itemize}
\setlength{\parsep}{0pt}
\setlength{\itemsep}{0pt}
\item Let $G=\U(m)\subset{\rm O}(2m)$. Then $G$ preserves a 2-form
$\om_0$ on $\R^{2m}$. If $g$ is a metric on $M$ with holonomy
$\U(m)$ then $g$ is {\it K\"ahler} with complex structure $J$, and 
the 2-form $\om$ on $M$ associated to $\om_0$ is the {\it K\"ahler 
form} of $g$. 

One can show that $\om$ is a calibration on $(M,g)$, and the 
calibrated submanifolds are exactly the {\it holomorphic curves} 
in $(M,J)$. More generally $\om^k/k!$ is a calibration on $M$ for 
$1\le k\le m$, and the corresponding calibrated submanifolds are 
the complex $k$-dimensional submanifolds of~$(M,J)$.
\item Let $G=\SU(m)\subset{\rm O}(2m)$. Then $G$ preserves a
{\it complex volume form} $\Om_0=\d z_1\w\cdots\w\d z_m$ on $\C^m$.
Thus a {\it Calabi--Yau $m$-fold\/} $(M,g)$ with $\Hol(g)=\SU(m)$
has a {\it holomorphic volume form} $\Om$. The real part $\Re\Om$
is a calibration on $M$, and the corresponding calibrated
submanifolds are called {\it special Lagrangian submanifolds}.
\item The group $G_2\subset{\rm O}(7)$ preserves a 3-form $\vp_0$ and 
a 4-form $*\vp_0$ on $\R^7$. Thus a Riemannian 7-manifold $(M,g)$ with 
holonomy $G_2$ comes with a 3-form $\vp$ and 4-form $*\vp$, which are
both calibrations. The corresponding calibrated submanifolds are
called {\it associative $3$-folds} and {\it coassociative $4$-folds}.
\item The group $\Spin(7)\subset{\rm O}(8)$ preserves a 4-form $\Om_0$
on $\R^8$. Thus a Riemannian 8-manifold $(M,g)$ with holonomy Spin(7) 
has a 4-form $\Om$, which is a calibration. We call $\Om$-submanifolds 
{\it Cayley $4$-folds}.
\end{itemize}

It is an important general principle that to each calibration
$\vp$ on an $n$-manifold $(M,g)$ with special holonomy we
construct in this way, there corresponds a constant calibration
$\vp_0$ on $\R^n$. Locally, $\vp$-submanifolds in $M$ will look 
very like $\vp_0$-submanifolds in $\R^n$, and have many of the
same properties. Thus, to understand the calibrated submanifolds 
in a manifold with special holonomy, it is often a good idea to 
start by studying the corresponding calibrated submanifolds of~$\R^n$. 

In particular, singularities of $\vp$-submanifolds in $M$ will be 
locally modelled on singularities of $\vp_0$-submanifolds in $\R^n$. 
(In the sense of Geometric Measure Theory, the {\it tangent cone}
at a singular point of a $\vp$-submanifold in $M$ is a conical
$\vp_0$-submanifold in $\R^n$.) So by studying singular
$\vp_0$-submanifolds in $\R^n$, we may understand the singular
behaviour of $\vp$-submanifolds in~$M$.

\subsection{Associative and coassociative submanifolds}
\label{ec53}

We now discuss the calibrated submanifolds of $G_2$-manifolds.

\begin{dfn} Let $(M,\vp,g)$ be a $G_2$-{\it manifold}, as
in \S\ref{ec22}. Then the 3-form $\vp$ is a {\it calibration}
on $(M,g)$. We define an {\it associative $3$-fold} in $M$ to
be a 3-submanifold of $M$ calibrated with respect to $\vp$.
Similarly, the Hodge star $*\vp$ of $\vp$ is a calibration
4-form on $(M,g)$. We define a {\it coassociative $4$-fold}
in $M$ to be a 4-submanifold of $M$ calibrated with respect
to~$*\vp$.
\label{ec5def2}
\end{dfn}

To understand these, it helps to begin with some calculations
on $\R^7$. Let the metric $g_0$, 3-form $\vp_0$ and 4-form
$*\vp_0$ on $\R^7$ be as in \S\ref{ec22}. Define an {\it
associative $3$-plane} to be an oriented 3-dimensional vector
subspace $V$ of $\R^7$ with $\vp_0\vert_V=\vol_V$, and a {\it
coassociative $4$-plane} to be an oriented 4-dimensional
vector subspace $V$ of $\R^7$ with $*\vp_0\vert_V=\vol_V$.
From \cite[Th.~IV.1.8, Def.~IV.1.15]{HaLa} we have:

\begin{prop} The family ${\mathcal F}^3$ of associative $3$-planes in
$\R^7$ and the family ${\mathcal F}^4$ of coassociative $4$-planes in
$\R^7$ are both isomorphic to $G_2/SO(4)$, with dimension~$8$.
\label{ec5prop2}
\end{prop}

Examples of an associative 3-plane $U$ and a coassociative
4-plane $V$ are
\e
U=\bigl\{(x_1,x_2,x_3,0,0,0,0):x_j\in\R\bigr\}
\;\>\text{and}\;\>
V=\bigl\{(0,0,0,x_4,x_5,x_6,x_7):x_j\in\R\bigr\}.
\label{ec5eq1}
\e
As $G_2$ acts {\it transitively} on the set of associative
3-planes by Proposition \ref{ec5prop2}, every associative
3-plane is of the form $\ga\cdot U$ for $\ga\in G_2$.
Similarly, every coassociative 4-plane is of the form
$\ga\cdot V$ for~$\ga\in G_2$.

Now $\vp_0\vert_V\equiv 0$. As $\vp_0$ is $G_2$-invariant,
this gives $\vp_0\vert_{\ga\cdot V}\equiv 0$ for all
$\ga\in G_2$, so $\vp_0$ restricts to zero on all
coassociative 4-planes. In fact the converse is true:
if $W$ is a 4-plane in $\R^7$ with $\vp_0\vert_W\equiv 0$,
then $W$ is coassociative with some orientation. From
this we deduce an alternative characterization of
coassociative 4-folds:

\begin{prop} Let\/ $(M,\vp,g)$ be a $G_2$-manifold, and\/ $L$
a $4$-dimensional submanifold of\/ $M$. Then $L$ admits an
orientation making it into a coassociative $4$-fold if and
only if\/~$\vp\vert_L\equiv 0$.
\label{ec5prop3}
\end{prop}

Trivially, $\vp\vert_L\equiv 0$ implies that $[\vp\vert_L]=0$
in $H^3(L,\R)$. Regard $L$ as an immersed 4-submanifold, with
immersion $\io:L\ra M$. Then $[\vp\vert_L]\in H^3(L,\R)$ is
unchanged under continuous variations of the immersion $\io$.
Thus, $[\vp\vert_L]=0$ is a necessary condition not just for
$L$ to be coassociative, but also for any isotopic 4-fold $N$
in $M$ to be coassociative. This gives a {\it topological
restriction} on coassociative 4-folds.

\begin{cor} Let\/ $(\vp,g)$ be a torsion-free $G_2$-structure on 
a $7$-manifold\/ $M$, and\/ $L$ a real\/ $4$-submanifold in $M$.
Then a necessary condition for $L$ to be isotopic to a coassociative
$4$-fold\/ $N$ in $M$ is that\/ $[\vp\vert_L]=0$ in~$H^3(L,\R)$.
\label{ec5cor}
\end{cor}

\subsection{Examples of associative 3-submanifolds}
\label{ec54}

Here are some sources of examples of associative 3-folds in~$\R^7$:
\begin{itemize}
\item Write $\R^7=\R\op\C^3$. Then $\R\t\Si$ is an associative
3-fold in $\R^7$ for any {\it holomorphic curve\/} $\Si$ in
$\C^3$. Also, if $L$ is any {\it special Lagrangian $3$-fold\/}
in $\C^3$ and $x\in\R$ then $\{x\}\t L$ is associative 3-fold
in $\R^7$. For examples of special Lagrangian 3-folds see
\cite[\S 9]{GHJ}, and references therein.
\item Bryant \cite[\S 4]{Brya1} studies compact Riemann surfaces
$\Si$ in ${\mathcal S}^6$ pseudoholomorphic with respect to the
almost complex structure $J$ on ${\mathcal S}^6$ induced by its
inclusion in $\Im{\mathbb O}\cong\R^7$. Then the cone on $\Si$
is an {\it associative cone\/} on $\R^7$. He shows that any
$\Si$ has a {\it torsion\/} $\tau$, a holomorphic analogue of
the Serret--Frenet torsion of real curves in~$\R^3$.

The torsion $\tau$ is a section of a holomorphic line bundle
on $\Si$, and $\tau=0$ if $\Si\cong\CP^1$. If $\tau=0$ then
$\Si$ is the projection to ${\mathcal S}^6=G_2/\SU(3)$ of
a {\it holomorphic curve\/} $\ti\Si$ in the {\it projective
complex manifold\/} $G_2/\U(2)$. This reduces the problem
of understanding null-torsion associative cones in $\R^7$ to
that of finding holomorphic curves $\ti\Si$ in $G_2/\U(2)$
satisfing a {\it horizontality condition}, which is a problem
in {\it complex algebraic geometry}. In integrable systems
language, null torsion curves are called {\it superminimal}.

Bryant also shows that {\it every\/} Riemann surface $\Si$ may
be embedded in ${\mathcal S}^6$ with null torsion in infinitely
many ways, of arbitrarily high degree. This shows that there
are {\it many associative cones in\/} $\R^7$, on {\it oriented
surfaces of every genus}. These provide many local models
for {\it singularities\/} of associative 3-folds.

Perhaps the simplest nontrivial example of a pseudoholomorphic
curve $\Si$ in ${\mathcal S}^6$ with null torsion is the
{\it Bor{\accent'27u}vka sphere} \cite{Boru}, which is an
${\mathcal S}^2$ orbit of an $\SO(3)$ subgroup of $G_2$
acting irreducibly on $\R^7$. Other examples are given by
Ejiri \cite[\S 5--\S 6]{Ejir}, who classifies pseudoholomorphic
${\mathcal S}^2$'s in ${\mathcal S}^6$ invariant under a $\U(1)$
subgroup of $G_2$, and Sekigawa~\cite{Seki}.
\item Bryant's paper is one of the first steps in the study of
associative cones in $\R^7$ using the theory of {\it integrable
systems}. Bolton et al.\ \cite{BVW}, \cite[\S 6]{BPW} use
integrable systems methods to prove important results on 
pseudoholomorphic curves $\Si$ in ${\mathcal S}^6$. When
$\Si$ is a torus $T^2$, they show it is of {\it finite type}
\cite[Cor.~6.4]{BPW}, and so can be classified in terms of
algebro-geometric {\it spectral data}, and perhaps even in
principle be written down explicitly.
\item {\it Curvature properties} of pseudoholomorphic curves
in ${\mathcal S}^6$ are studied by Hashimoto \cite{Hash} and
Sekigawa~\cite{Seki}.
\item Lotay \cite{Lota1} studies constructions for associative
3-folds $N$ in $\R^7$. These generally involve writing $N$ as
the total space of a 1-parameter family of surfaces $P_t$ in
$\R^7$ of a prescribed form, and reducing the condition for
$N$ to be associative to an o.d.e.\ in $t$, which can be
(partially) solved fairly explicitly.

Lotay also considers {\it ruled associative $3$-folds}
\cite[\S 6]{Lota1}, which are associative 3-folds $N$ in
$\R^7$ fibred by a 2-parameter family of affine straight
lines $\R$. He shows that any {\it associative cone} $N_0$
on a Riemann surface $\Si$ in ${\mathcal S}^6$ is the limit
of a 6-dimensional family of {\it Asymptotically Conical\/}
ruled associative 3-folds if $\Si\cong\CP^1$, and of a
2-dimensional family if~$\Si\cong T^2$.

Combined with the results of Bryant \cite[\S 4]{Brya1}
above, this yields many examples of generically nonsingular
Asymptotically Conical associative 3-folds in $\R^7$,
diffeomorphic to ${\mathcal S}^2\t\R$ or~$T^2\t\R$.
\end{itemize}

Examples of associative 3-folds in other explicit
$G_2$-manifolds, such as those of Bryant and Salamon
\cite{BrSa}, may also be constructed using similar
techniques. For finding associative 3-folds in
{\it nonexplicit\/} $G_2$-manifolds, such as the
compact examples of \S\ref{ec3} which are known
only through existence theorems, there is one
method \cite[\S 12.6]{Joyc5}, which we now explain.

Suppose $\ga\in G_2$ with $\ga^2=1$ but $\ga\ne 1$. Then
$\ga$ is conjugate in $G_2$ to
\begin{equation*}
(x_1,\dots,x_7)\longmapsto(x_1,x_2,x_3,-x_4,-x_5,-x_6,-x_7).
\end{equation*}
The fixed point set of this involution is the associative 3-plane
$U$ of \eq{ec5eq1}. It follows that any $\ga\in G_2$ with
$\ga^2=1$ but $\ga\ne 1$ has fixed point set an associative
3-plane. Thus we deduce~\cite[Prop.~10.8.1]{Joyc5}:

\begin{prop} Let\/ $(M,\vp,g)$ be a $G_2$-manifold, and\/
$\si:M\ra M$ be a nontrivial isometric involution with\/
$\si^*(\vp)=\vp$. Then $N=\bigl\{p\in M:\si(p)=p\bigr\}$
is an associative $3$-fold in~$M$.
\label{ec5prop4}
\end{prop}

Here a {\it nontrivial isometric involution} of $(M,g)$ is
a diffeomorphism $\si:M\ra M$ such that $\si^*(g)=g$, and
$\si\ne\id$ but $\si^2=\id$, where $\id$ is the identity on
$M$. Following \cite[Ex.~12.6.1]{Joyc5}, we can use the
proposition in to construct {\it examples} of compact
associative 3-folds in the compact 7-manifolds with
holonomy $G_2$ constructed in~\S\ref{ec3}.

\begin{ex} Let $T^7=\R^7/\Z^7$ and $\Ga$ be as in Example
\ref{ec3ex}. Define $\si:T^7\ra T^7$ by
\begin{equation*}
\si:(x_1,\ldots,x_7)\mapsto(x_1,x_2,x_3,\ha-x_4,-x_5,-x_6,-x_7).
\end{equation*}
Then $\si$ preserves $(\vp_0,g_0)$ and commutes with $\Ga$,
and so its action pushes down to $T^7/\Ga$. The fixed points
of $\si$ on $T^7$ are 16 copies of $T^3$, and $\si\de$ has no
fixed points in $T^7$ for all $\de\ne 1$ in $\Ga$. Thus the
fixed points of $\si$ in $T^7/\Ga$ are the image of the 16
$T^3$ fixed by $\si$ in~$T^7$.

But calculation shows that these 16 $T^3$ do not intersect the fixed 
points of $\al$, $\be$ or $\ga$, and that $\Ga$ acts freely on the 
set of 16 $T^3$ fixed by $\si$. So the image of the 16 $T^3$ in $T^7$ 
is 2 $T^3$ in $T^7/\Ga$, which do not intersect the singular set of 
$T^7/\Ga$, and which are {\it associative $3$-folds} in $T^7/\Ga$ by
Proposition~\ref{ec5prop4}.

Now the resolution of $T^7/\Ga$ to get a compact $G_2$-manifold
$(M,\ti\vp,\ti g)$ with $\Hol(\ti g)=G_2$ described in \S\ref{ec3}
may be done in a $\si$-{\it equivariant way}, so that $\si$ lifts
to $\si:M\ra M$ with $\si^*(\ti\vp)=\ti\vp$. The fixed points of
$\si$ in $M$ are again 2 copies of $T^3$, which are {\it
associative $3$-folds} by Proposition~\ref{ec5prop4}.
\label{ec5ex1}
\end{ex}

\subsection{Examples of coassociative 4-submanifolds}
\label{ec55}

Here are some sources of examples of coassociative 4-folds
in~$\R^7$:
\begin{itemize}
\item Write $\R^7=\R\op\C^3$. Then $\{x\}\t S$ is a
coassociative 4-fold in $\R^7$ for any {\it holomorphic
surface} $S$ in $\C^3$ and $x\in\R$. Also, $\R\t L$ is a
coassociative 4-fold in $\R^7$ for any {\it special
Lagrangian $3$-fold\/ $L$ in $\C^3$ with phase} $i$.
For examples of special Lagrangian 3-folds see
\cite[\S 9]{GHJ}, and references therein.
\item Harvey and Lawson \cite[\S IV.3]{HaLa} give examples
of coassociative 4-folds in $\R^7$ invariant under $\SU(2)$,
acting on $\R^7\cong\R^3\op\C^2$ as $\SO(3)=\SU(2)/\{\pm1\}$
on the $\R^3$ and $\SU(2)$ on the $\C^2$ factor. Such 4-folds
correspond to solutions of an o.d.e., which Harvey and
Lawson solve.
\item Mashimo \cite{Mash} classifies {\it coassociative cones}
$N$ in $\R^7$ with $N\cap{\mathcal S}^6$ homogeneous under a
3-dimensional simple subgroup $H$ of~$G_2$.
\item Lotay \cite{Lota2} studies 2-{\it ruled coassociative
$4$-folds} in $\R^7$, that is, coassociative 4-folds $N$ which
are fibred by a 2-dimensional family of affine 2-planes $\R^2$
in $\R^7$, with base space a Riemann surface $\Si$. He shows
that such 4-folds arise locally from data $\phi_1,\phi_2:\Si
\ra{\mathcal S}^6$ and $\psi:\Si\ra\R^7$ satisfying nonlinear
p.d.e.s similar to the Cauchy--Riemann equations.

For $\phi_1,\phi_2$ fixed, the remaining equations on $\psi$
are {\it linear}. This means that the family of 2-ruled
associative 4-folds $N$ in $\R^7$ asymptotic to a fixed
2-ruled coassociative cone $N_0$ has the structure of
a {\it vector space}. It can be used to generate families
of examples of coassociative 4-folds in~$\R^7$.
\end{itemize}

We can also use the fixed-point set technique of \S\ref{ec54} to
find examples of coassociative 4-folds in other $G_2$-manifolds.
If $\al:\R^7\ra\R^7$ is linear with $\al^2=1$ and $\al^*(\vp_0)=
-\vp_0$, then either $\al=-1$, or $\al$ is conjugate under an
element of $G_2$ to the map
\begin{equation*}
(x_1,\dots,x_7)\longmapsto(-x_1,-x_2,-x_3,x_4,x_5,x_6,x_7).
\end{equation*}
The fixed set of this map is the coassociative 4-plane $V$ 
of \eq{ec5eq1}. Thus, the fixed point set of $\al$ is 
either $\{0\}$, or a coassociative 4-plane in $\R^7$. So
we find~\cite[Prop.~10.8.5]{Joyc5}:

\begin{prop} Let\/ $(M,\vp,g)$ be a $G_2$-manifold, and\/
$\si:M\ra M$ an isometric involution with\/ $\si^*(\vp)=-\vp$.
Then each connected component of the fixed point set\/
$\bigl\{p\in M:\si(p)=p\bigr\}$ of\/ $\si$ is either a
coassociative $4$-fold or a single point.
\label{ec5prop5}
\end{prop}

Bryant \cite{Brya3} uses this idea to construct many
{\it local\/} examples of compact coassociative 4-folds
in $G_2$-manifolds.

\begin{thm}[Bryant \cite{Brya3}] Let\/ $(N,g)$ be a compact,
real analytic, oriented Riemannian $4$-manifold whose bundle
of self-dual\/ $2$-forms is trivial. Then $N$ may be embedded
isometrically as a coassociative $4$-fold in a $G_2$-manifold
$(M,\vp,g)$, as the fixed point set of an involution~$\si$. 
\label{ec5thm}
\end{thm}

Note here that $M$ need not be {\it compact}, nor $(M,g)$ {\it
complete}. Roughly speaking, Bryant's proof constructs $(\vp,g)$
as the sum of a power series on $\La^2_+T^*N$ converging
near the zero section $N\subset\La^2T^*N$, using the theory of
{\it exterior differential systems}. The involution $\si$ acts
as $-1$ on $\La^2_+T^*N$, fixing the zero section. One moral
of Theorem \ref{ec5thm} is that to be coassociative places
no significant local restrictions on a 4-manifold, other than
orientability.

Examples of {\it compact\/} coassociative 4-folds in {\it
compact\/} $G_2$-manifolds with holonomy $G_2$ are constructed
in \cite[\S 12.6]{Joyc5}, using Proposition \ref{ec5prop5}. Here
\cite[Ex.~12.6.4]{Joyc5} are examples in the $G_2$-manifolds
of~\S\ref{ec3}.

\begin{ex} Let $T^7=\R^7/\Z^7$ and $\Ga$ be as in Example
\ref{ec3ex}. Define $\si:T^7\ra T^7$ by
\begin{equation*}
\si:(x_1,\ldots,x_7)\mapsto(\ha-x_1,x_2,x_3,x_4,x_5,\ha-x_6,\ha-x_7).
\end{equation*}
Then $\si$ commutes with $\Ga$, preserves $g_0$ and takes $\vp_0$ to
$-\vp_0$. The fixed points of $\si$ in $T^7$ are 8 copies of $T^4$,
and the fixed points of $\si\al\be$ in $T^7$ are 128 points. If
$\de\in\Ga$ then $\si\de$ has no fixed points unless $\de=1,\al\be$.
Thus the fixed points of $\si$ in $T^7/\Ga$ are the image of the
fixed points of $\si$ and $\si\al\be$ in~$T^7$. 

Now $\Ga$ acts freely on the sets of 8 $\si$ $T^4$ and 128 $\si\al\be$
points. So the fixed point set of $\si$ in $T^7/\Ga$ is the union of 
$T^4$ and 16 isolated points, none of which intersect the singular set 
of $T^7/\Ga$. When we resolve $T^7/\Ga$ to get $(M,\ti\vp,\ti g)$ with
$\Hol(\ti g)=G_2$ in a $\si$-equivariant way, the action of $\si$ on
$M$ has $\si^*(\ti\vp)=-\ti\vp$, and again fixes $T^4$ and 16 points.
By Proposition \ref{ec5prop5}, this $T^4$ is {\it coassociative}.
\label{ec5ex2}
\end{ex}

More examples of associative and coassociative submanifolds
with different topologies are given in~\cite[\S 12.6]{Joyc5}.

\subsection{Cayley 4-folds}
\label{ec56}

The calibrated geometry of $\Spin(7)$ is similar to the $G_2$
case above, so we shall be brief.

\begin{dfn} Let $(M,\Om,g)$ be a $\Spin(7)$-{\it manifold}, as
in \S\ref{ec23}. Then the 4-form $\Om$ is a {\it calibration}
on $(M,g)$. We define a {\it Cayley $4$-fold} in $M$ to be a
4-submanifold of $M$ calibrated with respect to~$\Om$.
\label{ec5def3}
\end{dfn}

Let the metric $g_0$, and 4-form $\Om_0$ on $\R^8$ be as in
\S\ref{ec23}. Define a {\it Cayley $4$-plane} to be an oriented
4-dimensional vector subspace $V$ of $\R^8$ with $\Om_0\vert_V=
\vol_V$. Then we have an analogue of Proposition~\ref{ec5prop2}:

\begin{prop} The family $\mathcal F$ of Cayley $4$-planes in
$\R^8$ is isomorphic to $\Spin(7)/K$, where $K\cong
\bigl(\SU(2)\t\SU(2)\t\SU(2)\bigr)/\Z_2$ is a Lie subgroup 
of\/ $\Spin(7)$, and\/~$\dim{\mathcal F}=12$.
\label{ec5prop6}
\end{prop}

Here are some sources of examples of Cayley 4-folds in~$\R^8$:
\begin{itemize}
\item Write $\R^8=\C^4$. Then any {\it holomorphic surface} $S$
in $\C^4$ is Cayley in $\R^8$, and any {\it special Lagrangian
$4$-fold\/} $N$ in $\C^4$ is Cayley in~$\R^8$.

\noindent Write $\R^8=\R\t\R^7$. Then $\R\t L$ is Cayley for any
{\it associative $3$-fold\/} $L$ in~$\R^7$.
\item Lotay \cite{Lota2} studies 2-{\it ruled Cayley $4$-folds}
in $\R^8$, that is, Cayley 4-folds $N$ fibred by a 2-dimensional
family $\Si$ of affine 2-planes $\R^2$ in $\R^8$, as for the
coassociative case in \S\ref{ec55}. He constructs explicit
families of 2-ruled Cayley 4-folds in $\R^8$, including some
depending on an arbitrary holomorphic function $w:\C\ra\C$,
\cite[Th.~5.1]{Lota2}.
\end{itemize}

By the method of Propositions \ref{ec5prop4} and 
\ref{ec5prop5} one can prove~\cite[Prop.~10.8.6]{Joyc5}:

\begin{prop} Let\/ $(M,\Om,g)$ be a $\Spin(7)$-manifold,
and\/ $\si:M\ra M$ a nontrivial isometric involution with\/
$\si^*(\Om)=\Om$. Then each connected component of the fixed
point set\/ $\bigl\{p\in M:\si(p)=p\bigr\}$ is either a Cayley
$4$-fold or a single point.
\label{ec5prop7}
\end{prop}

Using this, \cite[\S 14.3]{Joyc5} constructs examples
of {\it compact\/} Cayley 4-folds in compact 8-manifolds
with holonomy~$\Spin(7)$.

\section{Deformations of calibrated submanifolds}
\label{ec6}

Finally we discuss {\it deformations} of associative,
coassociative and Cayley submanifolds. In \S\ref{ec61}
we consider the local equations for such submanifolds
in $\R^7$ and $\R^8$, following Harvey and Lawson
\cite[\S IV.2]{HaLa}. Then \S\ref{ec62} explains the
deformation theory of {\it compact\/} coassociative
4-folds, following McLean \cite[\S 4]{McLe}. This has
a particularly simple structure, as coassociative
4-folds are defined by the vanishing of $\vp$. The
deformation theory of compact associative 3-folds
and Cayley 4-folds is more complex, and is sketched
in~\S\ref{ec63}.

\subsection{Parameter counting and the local equations}
\label{ec61} 

We now study the local equations for 3- or 4-folds to be
(co)associative or Cayley.
\smallskip

\noindent{\it Associative 3-folds.}
The set of all 3-planes in $\R^7$ has dimension 12, and the set
of associative 3-planes in $\R^7$ has dimension 8 by Proposition
\ref{ec5prop2}. Thus the associative 3-planes are of {\it
codimension} 4 in the set of all 3-planes. Therefore the
condition for a 3-fold $L$ in $\R^7$ to be associative is 4
equations on each tangent space. The freedom to vary $L$ is
the sections of its normal bundle in $\R^7$, which is 4 real
functions. Thus, the deformation problem for associative
3-folds involves 4 {\it equations on $4$ functions}, so it
is a {\it determined\/} problem.

To illustrate this, let $f:\R^3\ra\H$ be a smooth function, written
\begin{equation*}
f(x_1,x_2,x_3)=f_0(x_1,x_2,x_3)+f_1(x_1,x_2,x_3)i+
f_2(x_1,x_2,x_3)j+f_3(x_1,x_2,x_3)k.
\end{equation*}
Define a 3-submanifold $L$ in $\R^7$ by
\begin{equation*}
L=\bigl\{\bigl(x_1,x_2,x_3,f_0(x_1,x_2,x_3),\ldots,
f_3(x_1,x_2,x_3)\bigr):x_j\in\R\bigr\}.
\end{equation*}
Then Harvey and Lawson \cite[\S IV.2.A]{HaLa} calculate the
conditions on $f$ for $L$ to be associative. With the
conventions of \S\ref{ec21}, the equation is
\e
i\frac{\pd f}{\pd x_1}+j\frac{\pd f}{\pd x_2}-k\frac{\pd f}{\pd x_3}
=C\Bigl(\frac{\pd f}{\pd x_1},\frac{\pd f}{\pd x_2},
\frac{\pd f}{\pd x_3}\Bigr),
\label{ec6eq1}
\e
where $C:\H\t\H\t\H\ra\H$ is a trilinear cross product.

Here \eq{ec6eq1} is 4 equations on 4 functions, as we claimed,
and is a {\it first order nonlinear elliptic p.d.e}. When
$f,\pd f$ are small, so that $L$ approximates the associative
3-plane $U$ of \eq{ec5eq1}, equation \eq{ec6eq1} reduces
approximately to the linear equation $i\frac{\pd f}{\pd x_1}+
j\frac{\pd f}{\pd x_2}-k\frac{\pd f}{\pd x_3}=0$, which is
equivalent to the {\it Dirac equation} on $\R^3$. More
generally, first order deformations of an associative
3-fold $L$ in a $G_2$-manifold $(M,\vp,g)$ correspond
to solutions of a {\it twisted Dirac equation} on~$L$.
\smallskip

\noindent{\it Coassociative 4-folds.}
The set of all 4-planes in $\R^7$ has dimension 12, and the set
of coassociative 4-planes in $\R^7$ has dimension 8 by Proposition
\ref{ec5prop2}. Thus the coassociative 4-planes are of {\it
codimension} 4 in the set of all 4-planes. Therefore the
condition for a 4-fold $N$ in $\R^7$ to be coassociative is 4
equations on each tangent space. The freedom to vary $N$ is
the sections of its normal bundle in $\R^7$, which is 3 real
functions. Thus, the deformation problem for coassociative
4-folds involves 4 {\it equations on $3$ functions}, so it
is an {\it overdetermined\/} problem.

To illustrate this, let $f:\H\ra\R^3$ be a smooth function, written
\begin{equation*}
f(x_0+x_1i+x_2j+x_3k)=(f_1,f_2,f_3)(x_0+x_1i+x_2j+x_3k).
\end{equation*}
Define a 4-submanifold $N$ in $\R^7$ by
\begin{equation*}
N=\bigl\{\bigl(f_1(x_0,\ldots,x_3),f_2(x_0,\ldots,x_3),
f_3(x_0,\ldots,x_3),x_0,\ldots,x_3\bigr):x_j\in\R\bigr\}.
\end{equation*}
Then Harvey and Lawson \cite[\S IV.2.B]{HaLa} calculate
the conditions on $f$ for $N$ to be coassociative. With
the conventions of \S\ref{ec21}, the equation is
\e
i\pd f_1+j\pd f_2-k\pd f_3=C(\pd f_1,\pd f_2,\pd f_3),
\label{ec6eq2}
\e
where the derivatives $\pd f_j=\pd f_j(x_0+x_1i+x_2j+x_3k)$
are interpreted as functions $\H\ra\H$, and $C$ is as
in \eq{ec6eq1}. Here \eq{ec6eq2} is 4 equations on 3
functions, as we claimed, and is a {\it first order
nonlinear overdetermined elliptic p.d.e}.
\smallskip

\noindent{\it Cayley 4-folds.}
The set of all 4-planes in $\R^8$ has dimension 16, and the
set of Cayley 4-planes in $\R^8$ has dimension 12 by
Proposition \ref{ec5prop6}, so the Cayley 4-planes are of
{\it codimension} 4 in the set of all 4-planes. Therefore
the condition for a 4-fold $K$ in $\R^8$ to be Cayley is
4 equations on each tangent space. The freedom to vary $K$
is the sections of its normal bundle in $\R^8$, which is
4 real functions. Thus, the deformation problem for
Cayley 4-folds involves 4 {\it equations on $4$ functions},
so it is a {\it determined\/} problem.

Let $f=f_0+f_1i+f_2j+f_3k=f(x_0+x_1i+x_2j+x_3k):\H\ra\H$ be
smooth. Choosing signs for compatibility with \eq{ec2eq2},
define a 4-submanifold $K$ in $\R^8$ by
\begin{align*}
K=\bigl\{\bigl(-x_0,x_1,&x_2,x_3,f_0(x_0+x_1i+x_2j+x_3k),
-f_1(x_0+x_1i+x_2j+x_3k),\\
-&f_2(x_0+x_1i+x_2j+x_3k),f_3(x_0+x_1i+x_2j+x_3k)\bigr):x_j\in\R\bigr\}.
\end{align*}
Following
\cite[\S IV.2.C]{HaLa}, the equation for $K$ to be Cayley is
\e
\frac{\pd f}{\pd x_0}+i\frac{\pd f}{\pd x_1}+j\frac{\pd f}{\pd x_2}
+k\frac{\pd f}{\pd x_3}=C(\pd f),
\label{ec6eq3}
\e
for $C:\H\ot_{\R}\H\ra\H$ a homogeneous cubic polynomial. This is
4 equations on 4 functions, as we claimed, and is a first-order
nonlinear elliptic p.d.e.\ on $f$. The linearization at $f=0$
is equivalent to the {\it positive Dirac equation} on $\R^4$.
More generally, first order deformations of a Cayley 4-fold
$K$ in a $\Spin(7)$-manifold $(M,\Om,g)$ correspond to
solutions of a {\it twisted positive Dirac equation} on~$K$.

\subsection{Deformation theory of coassociative 4-folds}
\label{ec62}

Here is the main result in the deformation theory of coassociative
4-folds, proved by McLean \cite[Th.~4.5]{McLe}. As our sign
conventions for $\vp_0,*\vp_0$ in \eq{ec2eq1} are different
to McLean's, we use self-dual 2-forms in place of McLean's
anti-self-dual 2-forms.

\begin{thm} Let\/ $(M,\vp,g)$ be a $G_2$-manifold, and\/ $N$ a
compact coassociative $4$-fold in $M$. Then the moduli space
${\mathcal M}_N$ of coassociative $4$-folds isotopic to $N$ in $M$
is a smooth manifold of dimension~$b^2_+(N)$.
\label{ec6thm1}
\end{thm}

\begin{proof}[Sketch proof] Suppose for simplicity that $N$ is an
embedded submanifold. There is a natural orthogonal decomposition 
$TM\vert_N=TN\op\nu$, where $\nu\ra N$ is the {\it normal bundle} 
of $N$ in $M$. There is a natural isomorphism $\nu\cong\La^2_+T^*N$,
constructed as follows. Let $x\in N$ and $V\in\nu_x$. Then
$V\in T_xM$, so $V\cdot\vp\vert_x\in\La^2T_x^*M$, and
$(V\cdot\vp\vert_x)\vert_{T_xN}\in\La^2T_x^*N$. It turns out
that $(V\cdot\vp\vert_x)\vert_{T_xN}$ actually lies in
$\La^2_+T_x^*N$, the bundle of {\it self-dual\/ $2$-forms} on
$N$, and that the map $V\mapsto(V\cdot\vp\vert_x)\vert_{T_xN}$
defines an {\it isomorphism}~$\nu\,\smash{{\buildrel\cong\over\longra}}
\,\La^2_+T^*N$.

Let $T$ be a small {\it tubular neighbourhood} of $N$ in $M$. Then 
we can identify $T$ with a neighbourhood of the zero section in
$\nu$, using the exponential map. The isomorphism $\nu\cong\La^2_+
T^*N$ then identifies $T$ with a neighbourhood $U$ of the zero
section in $\La^2_+T^*N$. Let $\pi:T\ra N$ be the obvious projection.

Under this identification, submanifolds $N'$ in $T\subset M$
which are $C^1$ close to $N$ are identified with the {\it graphs}
$\Ga(\al)$ of small smooth sections $\al$ of $\La^2_+T^*N$ lying
in $U$. Write $C^\iy(U)$ for the subset of the vector space of
smooth self-dual 2-forms $C^\iy(\La^2_+T^*N)$ on $N$ lying in
$U\subset\La^2_+T^*N$. Then for each $\al\in C^\iy(U)$ the graph
$\Ga(\al)$ is a 4-submanifold of $U$, and so is identified with
a 4-submanifold of $T$. We need to know: which 2-forms $\al$ 
correspond to {\it coassociative} 4-folds $\Ga(\al)$ in~$T$?

Well, $N'$ is coassociative if $\vp\vert_{N'}\equiv 0$. Now
$\pi\vert_{N'}:N'\ra N$ is a diffeomorphism, so we can push
$\vp\vert_{N'}$ down to $N$, and regard it as a function
of $\al$. That is, we define
\e
P:C^\iy(U)\longra C^\iy(\La^3T^*N) \quad\text{by}\quad
P(\al)=\pi_*(\vp\vert_{\Ga(\al)}).
\label{ec6eq4}
\e
Then the moduli space ${\mathcal M}_N$ is locally isomorphic near
$N$ to the set of small self-dual 2-forms $\al$ on $N$ with
$\vp\vert_{\Ga(\al)}\equiv 0$, that is, to a neighbourhood of
0 in~$P^{-1}(0)$.

To understand the equation $P(\al)=0$, note that at
$x\in N$, $P(\al)\vert_x$ depends on the tangent space
to $\Ga(\al)$ at $\al\vert_x$, and so on $\al\vert_x$
and $\nabla\al\vert_x$. Thus the functional form of $P$ is
\begin{equation*}
P(\al)\vert_x=F\bigl(x,\al\vert_x,\nabla\al\vert_x\bigr)
\quad\text{for $x\in N$,}
\end{equation*}
where $F$ is a smooth function of its arguments. Hence $P(\al)=0$
is a {\it nonlinear first order p.d.e.} in $\al$. The {\it
linearization} $\d P(0)$ of $P$ at $\al=0$ turns out to be
\begin{equation*}
\d P(0)(\be)=\lim_{\ep\ra 0}\bigl(\ep^{-1}P(\ep\be)\bigr)=\d\be.
\end{equation*}

Therefore $\Ker(\d P(0))$ is the vector space ${\mathcal H}^2_+$
of {\it closed self-dual\/ $2$-forms} $\be$ on $N$, which by
Hodge theory is a finite-dimensional vector space isomorphic
to $H^2_+(N,\R)$, with dimension $b^2_+(N)$. This is the
{\it Zariski tangent space} of ${\mathcal M}_N$ at $N$, the
{\it infinitesimal deformation space} of $N$ as a
coassociative 4-fold.

To complete the proof we must show that ${\mathcal M}_N$ is
locally isomorphic to its Zariski tangent space ${\mathcal H}^2_+$,
and so is a smooth manifold of dimension $b^2_+(N)$. To do this
rigorously requires some technical analytic machinery, which
is passed over in a few lines in \cite[p.~731]{McLe}. Here is
one way to do it. 

Because $C^\iy(\La^2_+T^*N),C^\iy(\La^3T^*N)$ are not
{\it Banach spaces}, we extend $P$ in \eq{ec6eq4} to act
on {\it H\"older spaces} $C^{k+1,\ga}(\La^2_+T^*N),C^{k,\ga}
(\La^3T^*N)$ for $k\ge 1$ and $\ga\in(0,1)$, giving
\begin{equation*}
P_{k,\ga}:C^{k+1,\ga}(U)\longra C^{k,\ga}(\La^3T^*N)
\quad\text{defined by}\quad
P_{k,\ga}(\al)=\pi_*(\vp\vert_{\Ga(\al)}).
\end{equation*}
Then $P_{k,\ga}$ is a smooth map of Banach manifolds. Let
$V_{k,\ga}\subset C^{k,\ga}(\La^3T^*N)$ be the Banach subspace
of {\it exact\/} $C^{k,\ga}$ 3-forms on~$N$.

As $\vp$ is closed, $\vp\vert_N\equiv 0$, and $\Ga(\al)$ is
isotopic to $N$, we see that $\vp\vert_{\Ga(\al)}$ is an {\it
exact\/} 3-form on $\Ga(\al)$, so that $P_{k,\ga}$ maps into
$V_{k,\ga}$. The linearization
\begin{equation*}
\d P_{k,\ga}(0):C^{k+1,\ga}(\La^2_+T^*N)\longra V_{k,\ga},
\qquad \d P_{k,\ga}(0):\be\longmapsto\d\be
\end{equation*}
is then {\it surjective} as a map of Banach spaces. (To prove
this requires a discursion, using elliptic regularity results
for~$\d+\d^*$.)

Thus, $P_{k,\ga}:C^{k+1,\ga}(U)\ra V_{k,\ga}$ is a smooth
map of Banach manifolds, with $\d P_{k,\ga}(0)$ surjective. The
{\it Implicit Function Theorem for Banach spaces} now implies
that $P_{k,\ga}^{-1}(0)$ is near 0 a smooth submanifold of
$C^{k+1,\ga}(U)$, locally isomorphic to $\Ker(\d P_{k,\ga}(0))$.
But $P_{k,\ga}(\al)=0$ is an {\it overdetermined elliptic
equation} for small $\al$, and so elliptic regularity implies
that solutions $\al$ are smooth. Therefore $P_{k,\ga}^{-1}(0)=
P^{-1}(0)$ near 0, and similarly $\Ker(\d P_{k,\ga}(0))=
\Ker(\d P(0))={\mathcal H}^2_+$. This completes the proof.
\end{proof}

Here are some remarks on Theorem~\ref{ec6thm1}.
\begin{itemize}
\item This proof relies heavily on Proposition \ref{ec5prop3},
that a 4-fold $N$ in $M$ is coassociative if and only if
$\vp\vert_N\equiv 0$, for $\vp$ a closed 3-form on $M$. The
consequence of this is that the deformation theory of compact
coassociative 4-folds is {\it unobstructed}, and the
moduli space is {\it always} a smooth manifold with
dimension given by a topological formula.

{\it Special Lagrangian $m$-folds} of Calabi-Yau $m$-folds
can also be defined in terms of the vanishing of closed forms,
and their deformation theory is also unobstructed, as in
\cite[\S 3]{McLe} and \cite[\S 10.2]{GHJ}. However, associative
3-folds and Cayley 4-folds cannot be defined by the vanishing
of closed forms, and we will see in \S\ref{ec63} that this
gives their deformation theory a different flavour.
\item We showed in \S\ref{ec61} that the condition for a
4-fold $N$ in $M$ to be coassociative is locally 4 equations
on 3 functions, and so is {\it overdetermined}. However,
Theorem \ref{ec6thm1} shows that coassociative 4-folds
have {\it unobstructed\/} deformation theory, and often
form {\it positive-dimensional\/} moduli spaces. This
seems very surprising for an overdetermined equation.

The explanation is that the condition $\d\vp=0$ acts
as an {\it integrability condition} for the existence
of coassociative 4-folds. That is, since closed 3-forms on
$N$ essentially depend locally only on 3 real parameters, not
4, as $\vp$ is closed the equation $\vp\vert_N\equiv 0$ is
in effect only 3 equations on $N$ rather than 4, so we can
think of the deformation theory as really controlled by a
determined elliptic equation.

Therefore $\d\vp=0$ is essential for Theorem \ref{ec6thm1}
to work. In `almost $G_2$-manifolds' $(M,\vp,g)$ with
$\d\vp\ne 0$, the deformation problem for coassociative
4-folds is overdetermined and obstructed, and generically
there would be no coassociative 4-folds.
\item In Example \ref{ec5ex2} we constructed an example of a
compact coassociative 4-fold $N$ diffeomorphic to $T^4$ in a
compact $G_2$-manifold $(M,\vp,g)$. By Theorem \ref{ec6thm1},
$N$ lies in a {\it smooth\/ $3$-dimensional family} of
coassociative $T^4$'s in $M$. Locally, these may form a
{\it coassociative fibration} of~$M$.
\end{itemize}

Now suppose $\bigl\{(M,\vp_t,g_t):t\in(-\ep,\ep)\bigr\}$ is a
smooth 1-parameter family of $G_2$-manifolds, and $N_0$ a
compact coassociative 4-fold in $(M,\vp_0,g_0)$. When can
we extend $N_0$ to a smooth family of coassociative 4-folds
$N_t$ in $(M,\vp_t,g_t)$ for small $t$? By Corollary
\ref{ec5cor}, a necessary condition is that $[\vp_t
\vert_{N_0}]=0$ for all $t$. Our next result shows that
locally, this is also a {\it sufficient\/} condition.
It can be proved using similar techniques to Theorem 
\ref{ec6thm1}, though McLean did not prove~it.

\begin{thm} Let\/ $\bigl\{(M,\vp_t,g_t):t\in(-\ep,\ep)\bigr\}$
be a smooth\/ $1$-parameter family of\/ $G_2$-manifolds, and\/
$N_0$ a compact coassociative $4$-fold in $(M,\vp_0,g_0)$.
Suppose that\/ $[\vp_t\vert_{N_0}]=0$ in $H^3(N_0,\R)$ for
all\/ $t\in(-\ep,\ep)$. Then $N_0$ extends to a smooth\/
$1$-parameter family $\bigl\{N_t:t\in(-\de,\de)\bigr\}$,
where $0<\de\le\ep$ and\/ $N_t$ is a compact coassociative
$4$-fold in~$(M,\vp_t,g_t)$.
\label{ec6thm2}
\end{thm}

\subsection{Deformations of associative 3-folds and Cayley 4-folds}
\label{ec63}

Associative 3-folds and Cayley 4-folds cannot be defined in terms
of the vanishing of closed forms, and this gives their deformation
theory a different character to the coassociative case. Here is how
the theories work, drawn mostly from McLean~\cite[\S 5--\S 6]{McLe}.

Let $N$ be a compact associative 3-fold or Cayley 4-fold in a 
7- or 8-manifold $M$. Then there are vector bundles $E,F\ra N$ 
with $E\cong\nu$, the normal bundle of $N$ in $M$, and a 
first-order elliptic operator $D_N:C^\iy(E)\ra C^\iy(F)$ on $N$. 
The {\it kernel\/} $\Ker D_N$ is the set of {\it infinitesimal 
deformations} of $N$ as an associative 3-fold or Cayley 4-fold. 
The {\it cokernel\/} $\Coker D_N$ is the {\it obstruction space}
for these deformations. 

Both are finite-dimensional vector spaces, and
\begin{equation*}
\dim\Ker D_N-\dim\Coker D_N=\ind(D_N),
\end{equation*}
the {\it index} of $D_N$. It is a topological invariant, given in
terms of characteristic classes by the {\it Atiyah--Singer Index 
Theorem}. In the associative case we have $E\cong F$, and $D_N$ is
anti-self-adjoint, so that $\Ker(D_N)\cong\Coker(D_N)$ and
$\ind(D_N)=0$ automatically. In the Cayley case we have
\begin{equation*}
\ind(D_N)=\tau(N)-\ha\chi(N)-\ha [N]\cdot[N],
\end{equation*}
where $\tau$ is the signature, $\chi$ the Euler characteristic
and $[N]\cdot[N]$ the self-intersection of~$N$.

In a {\it generic} situation we expect $\Coker D_N=0$, and then 
deformations of $N$ will be unobstructed, so that the moduli space 
${\mathcal M}_N$ of associative or Cayley deformations of $N$ will 
locally be a smooth manifold of dimension $\ind(D_N)$. However,
in nongeneric situations the obstruction space may be nonzero,
and then the moduli space may not be smooth, or may have a
larger than expected dimension.

This general structure is found in the deformation theory of
other important mathematical objects --- for instance,
pseudo-holomorphic curves in almost complex manifolds, and
instantons and Seiberg--Witten solutions on 4-manifolds.
In each case, the moduli space is only smooth with
topologically determined dimension under a {\it genericity
assumption} which forces the obstructions to vanish.

\end{document}